\documentclass[10pt,a4paper]{article}
\usepackage{CJK}
\usepackage{mathrsfs}
\usepackage{amsfonts}
\begin{document}

\begin{center}
 {\huge\textbf{A New Operator Theory of Linear Partial Differential Equations}}
\vspace{0.3cm}

\footnotetext{\hspace*{-.45cm}\footnotesize $^\dag$ Corresponding
author E-mail: yananbiguangqing@sohu.com
\\${}^{\rm b)}$ E-mail: yuekaifly@163.com}

\begin{center}
\rm Guang-Qing Bi $^{\rm a)\dagger}$, \ \ Yue-Kai Bi $^{\rm b)}$
\end{center}

\end{center}

\noindent
\begin{abstract} We first strictly expressed the basic notions and research methods of abstract operators, which systematically expounded the main results of abstract operator theory. By combining abstract operators with the Laplace transform, we can easily apply this Laplace transform to $n+1$ dimensional partial differential equations. Further, all the analytic solutions to an initial value problem of an arbitrary order linear partial differential equation are expressed in these abstract operators. By writing abstract operators in this class into integral forms, the solutions in operator form are represented into integral forms. We thus solved the important problem of representing the solutions of linear higher-order partial differential equations into the integrations of some given functions. By introduction of abstract operators on Hilbert space, we further discuss the solvability of initial-boundary value problem for the linear higher-order partial differential equations.
\end{abstract}


\noindent{\bf MSC(2010) Subject Classification:}\ \ 47F05, 35A10, 35G16, 44A10

\section{Introduction}
\noindent

In 1960s, the general theory of linear partial differential
equations made important progress by using the generalized function
and its Fourier transform, further, the theory of
pseudo-differential operators. The concept of pseudo-differential operators is defined as (See \cite{chen})

If $a(x,\xi)\in{C}^\infty(\mathbb{R}_x^n\times\mathbb{R}_\xi^n)$,
and for arbitrary $\alpha,\beta\in\mathbb{N}^n$ and an real number
$m$,
\[|\partial_\xi^\alpha\partial_x^\beta{a}(x,\xi)|\leq{C}_{\alpha,\beta}(1+|\xi|)^{m-|\alpha|}\]
is tenable, where $C_{\alpha,\beta}$ is a constant, then the linear
continuous mapping $A$ :
$\mathscr{S}(\mathbb{R}^n)\rightarrow\mathscr{S}(\mathbb{R}^n)$ can
be defined as
\begin{equation}\label{chen1}
Au(x)=(2\pi)^{-n}\int{e}^{i\xi{x}}a(x,\xi)\hat{u}(\xi)d\xi,
\end{equation}
which is called the pseudo-differential operators, denoted by
$a(x,D)$, where $a(x,\xi)$ is the symbols of $a(x,D)$, $D$ denotes
\[D^\alpha=D_1^{\alpha_1}\cdots{D}_n^{\alpha_n}=\left(\frac{1}{i}\partial_{x_1}\right)^{\alpha_1}\left(\frac{1}{i}\partial_{x_2}\right)^{\alpha_2}\cdots
\left(\frac{1}{i}\partial_{x_n}\right)^{\alpha_n}.\]

This definition means the following equality is tenable,
namely
\begin{equation}\label{chen2}
a(x,D)e^{i\xi{x}}=a(x,\xi)e^{i\xi{x}},
\end{equation}
where $\xi{x}=\langle\xi,x\rangle=\xi_1x_1+\xi_2x_2+\cdots+\xi_nx_n$. Obviously, the properties of $a(x,D)$ depend on $a(x,\xi)$. equality
(\ref{chen2}) indicates the mapping relation between $a(x,D)$ and
$a(x,\xi)$. For a given function $a(x,\xi)$, if the
corresponding algorithms of $a(x,D)$ can be
derived from this type of mapping relations, to determine the domain
and range of $a(x,D)$, then this type of mapping relations
can be the best definition of pseudo-differential operators.
Luckily, this idea totally works by using the Analytic
continuous fundamental theorem. In 1997, the author delivered the
Analytic continuous fundamental theorem, based on which the author
introduced the concept of abstract operators, and derived the
algorithms of five types of abstract operators such as (See \cite{bi97})
\[\exp(h\partial_x),\;\sin(h\partial_x),\;\cos(h\partial_x),\;\sinh(h\partial_x)\;\mbox{and}\;\cosh(h\partial_x).\]
Where $h\partial_x=\langle{h,\partial_x}\rangle=h_1\partial_{x_1}+h_2\partial_{x_2}+\cdots+h_n\partial_{x_n},\;\cos(ih\partial_x)=\cosh(h\partial_x)$,
$\sin(ih\partial_x)=i\sinh(h\partial_x)$. Further, for more complex abstract operators such as
\[\exp(tP(\partial_x)),\;\cos(tP(\partial_x)^{1/2}),\;\frac{\sin(tP(\partial_x)^{1/2})}{P(\partial_x)^{1/2}},\;
\cosh(tP(\partial_x)^{1/2})\;\mbox{and}\;\frac{\sinh(tP(\partial_x)^{1/2})}{P(\partial_x)^{1/2}},\]
we can easily establish their algorithms (See \cite{bi97}, \cite{bi99}), thus the general solving
procedure of initial value problem for linear partial differential
equations is derived clearly (See \cite{bi97}-\cite{bi11}).
The abstract operators mainly aim at promoting partial differential operators in a different way than pseudo-differential operators to establish theories of partial differential equations.

Without doubt, the theory of abstract operators is a new and promising research field. The theory and methodology of abstract operators are gradually integrating with other mathematical branches. We believe the position and significance of abstract operators in the development of mathematics science will be generally acknowledged.

\section{Theory of abstract operators}
\noindent

The mathematical concept of abstract operators has seldom been understood since it in a minimalist form first introduced in 1997 (See \cite{bi97}). Therefore our primary task is to systematically illustrate fundamental theories of abstract operators.
\subsection{Basic notions}
\noindent

\textbf{Definition 2.1.} Within each convergence
circle of analytic functions, if the effects on the series term by
term from the linear operator converge uniformly to the effects on
the sum function, then the operator is called having the analytic
continuity.

\textbf{Definition 2.2.} $x^\alpha,\,x\in{\mathbb{R}^n}$ is called the base function of analytic function, its exponential
form is $e^{\xi{x}}=e^{\langle\xi,x\rangle}$. Where $\alpha\in\mathbb{N}^n,
\;\xi\in{\mathbb{R}_n}$ are respectively called
the character of $x^\alpha$ and $e^{\xi{x}}$.

\textbf{Definition 2.3.} A class of linear operators with the analytic continuity is
called the abstract operators, denoted by
$f(t,x\partial_x)$,
$x\in\Omega,\,t\in(0,T]$, which is also the linear continuous mapping $f(t,x\partial_x)$ :
$C^\infty((0,T]\times\Omega)\rightarrow{C^\infty((0,T]\times\Omega)}$, and if it acts on the base function
$x^\alpha$, we have
\[f(t,\,x\partial_x)x^\alpha=f(t,\alpha)x^\alpha,\quad\forall{f(t,\alpha)}\in{C^\infty}((0,T]\times\mathbb{R}_n).\]
Which is called the abstract operators taking
$x\partial_x=(x_1\partial_{x_1},\cdots,x_n\partial_{x_n})=(x_1\partial_1,\cdots,x_n\partial_n)$ as the operator element. Where
$f(t,\alpha),\alpha\in\mathbb{N}^n$ is called the symbols of the abstract operators
$f(t,x\partial_x)$.

\textbf{Definition 2.4.} Let $T>0,\,\Omega\subset\mathbb{R}^n$ be an open
set, a class of linear operators is called the abstract operators taking $\partial_x$ as the operator element, denoted by
$f(t,\partial_x),\,x\in\Omega,\,t\in(0,T]$, which is also the linear continuous mapping $f(t,\partial_x)$ :
$C^\infty((0,T]\times\Omega)\rightarrow{C^\infty((0,T]\times\Omega)}$, and if it
acts on the exponential base function $e^{\xi{x}}$, we have
\[f(t,\partial_x)e^{\xi{x}}=f(t,\xi)e^{\xi{x}},\quad\forall{f(t,\xi)}\in{C^\infty}((0,T]\times\mathbb{R}_n).\]
Where $f(t,\xi),\,\xi\in\mathbb{R}_n$ is called the symbols of the abstract operators $f(t,\partial_x)$.

\textbf{Analytic continuous fundamental theorem.} Let $A,B$ be the
linear operators with the analytic continuity, if there are
functions $y(x)\in\mathbb{R}^n,\;x\in\mathbb{R}^n$, making one of
the following two operator equalities
\[Ax^\alpha=By^\alpha\quad\mbox{or}\quad Ae^{\xi{x}}=Be^{\xi{y}}\]
 tenable, and the expressions of $A,B$ do not explicitly contain the
character $\alpha\in\mathbb{N}^n,\,\xi\in\mathbb{R}_n$, then
\[Af(x)=Bf(y),\qquad\forall f(x),f(y)\in\mathscr{D}(A)\cap\mathscr{D}(B),\]
where $\mathscr{D}(A)$ and $\mathscr{D}(B)$ are domain of linear operators. Especially, the Analytic continuous fundamental theorem has been
proved by Guang-Qing Bi in 1997 (See \cite{bi97}).

\textbf{Definition 2.5.} A linear operator equality, if it's true for
arbitrary functions in a certain range, then
it is called an operator formula. Particularly, the operator
formulas that determine the domain and range of abstract operators
are called the algorithms of the abstract operators.

\textbf{Definition 2.6.} The relational expression between each
component of the character
$\alpha=(\alpha_1,\alpha_2,\ldots,\alpha_n)$ or
$\xi=(\xi_1,\xi_2,\ldots,\xi_n)$ is called the characteristic
equation. In general, an operator formula becomes the corresponding
characteristic equation, when where arbitrary function of which is
the base function.

As the base function can be expressed by both
power function and exponential function, a single operator formula
can be corresponding to two different characteristic equations,
conversely, a single characteristic equation can also be
corresponding to two different operator formulas.

Constructing the operator equality tenable for the base function by the characteristic equation through the definition of abstract operators, and then deducing that it is also tenable for arbitrary analytic functions according to the Analytic continuous fundamental theorem, thus we derive new operator formula. Finding or establishing a new operator formula requires the knowledge of the corresponding characteristic equation in advance, without knowing the specific form of the new operator formula. Therefore, it all boils down to seek or construct appropriate characteristic equations. The key to transform the characteristic equation to the corresponding operator formula, is constructing the operator equality true for the base function by using the specific form of the characteristic equation and the definition of abstract operators, and also, only when the operators constructed are all linear and the expressions of the linear operators do not explicitly contain the character of the base function can we derive that the operator equality is not only true for the base function, but also for arbitrary analytic functions, according to the Analytic continuous fundamental theorem.

\textbf{Definition 2.7.} Let $A$ be a linear operator having the
analytic continuity, if there is another
linear operator, denoted by $A^{-1}$, making $AA^{-1}=A^{-1}A=I$,
then $A^{-1}$ is called the inverse operator of $A$.

By Definition 2.4 and the Analytic continuous fundamental theorem we obtain:

\textbf{Corollary 2.1.} The operator algebras formed
by all abstract operators $f(t,\partial_x)$, are isomorphic to
the algebras formed by their symbols $f(t,\xi)$. This isomorphism is
determined by $f(t,\partial_x)\leftrightarrow{f}(t,\xi)$, and
\[f(t,\partial_x)\pm{g}(t,\partial_x)\leftrightarrow{f}(t,\xi)\pm{g}(t,\xi),\]
$$f(t,\partial_x)\circ{g}(t,\partial_x)\leftrightarrow{f}(t,\xi)g(t,\xi).$$
Especially, the abstract operators $f(\partial_x)$ and
$g(\partial_x),\;x\in\mathbb{R}^n$ are each other's inverse
operators, if and only if their symbols $f(\xi)$ and $g(\xi)$ satisfy
$f(\xi)g(\xi)=1,\;\xi\in\mathbb{R}_n$.

For instance, $\cos(ix\partial_y)=\cosh(x\partial_y),\;\sin(ix\partial_y)=i\sinh(x\partial_y)$.
$\partial_{x_i}\sin(x\partial_y)=\partial_{y_i}\cos(x\partial_y)$,
$\partial_{x_i}\cos(x\partial_y)=-\partial_{y_i}\sin(x\partial_y)$,
where $x\partial_y=\langle{x,\partial_y}\rangle=x_1\partial_{y_1}+\cdots+x_n\partial_{y_n}$.

By Definition 2.4 we get:

\textbf{Corollary 2.2.} Let
$g(x)\in{C^\infty(\mathbb{R}^n)}$, if
$g(\partial_\xi)(e^{\xi{x}}f(\xi)),\,\xi\in\mathbb{R}_n$ is continuous at $\xi=\xi_0$, then
\begin{equation}\label{e0}
f(\partial_x)(e^{\xi_0x}g(x))=g(\partial_\xi)(e^{\xi{x}}f(\xi))|_{\xi=\xi_0}.
\end{equation}

By Definition 2.4 and the Analytic continuous fundamental theorem we get:

\textbf{Example 2.1.} Let $x\in\mathbb{R}^n,\;\xi\in\mathbb{R}_n$, if
$f(\xi+\partial_x)g(x)$ is continuous, then
\begin{equation}\label{e1}
f(\partial_x)(e^{\xi{x}}g(x))=e^{\xi{x}}f(\xi+\partial_x)g(x).
\end{equation}

Similarly, if $e^{-x\partial_y}(f(\partial_x)e^{x\partial_y}g(x,y))$ is continuous, then
\begin{equation}\label{e3}
f(\partial_x+\partial_y)g(x,y)=e^{-x\partial_y}(f(\partial_x)e^{x\partial_y}g(x,y)),
\end{equation}
where
$x\partial_y=x_1\partial_{y_1}+x_2\partial_{y_2}+\cdots+x_n\partial_{y_n},\;e^{x\partial_y}g(x,y)=g(x,y+x)$.

\textbf{Example 2.2.} Let $f(x),x\in\mathbb{R}^n,\,\lambda>0$,
then we have
\begin{equation}\label{e2}
 e^{\lambda|\partial_x|^2}f(x)=\frac{1}{2^n(\lambda\pi)^{n/2}}\int_{\mathbb{R}^n}f(\eta)\exp\left(-\frac{|\eta-x|^2}{4\lambda}\right)d\eta,
 \quad\lambda>0.
\end{equation}

\textbf{Proof.} If taking the following integral formula
\[e^{\lambda\xi_i^2}=\frac{1}{2\sqrt{\pi}}\int^\infty_{-\infty}e^{-\zeta^2/4+\sqrt{\lambda}\,\xi_i\zeta}d\zeta,\;\forall\lambda>0\]
as its characteristic equation, according to the Analytic continuous fundamental theorem we have
\begin{eqnarray*}
  \exp\left(\lambda\frac{\partial^2}{\partial{x_i^2}}\right)f(x)&=&
  \frac{1}{2\sqrt{\pi}}\int^\infty_{-\infty}e^{-\zeta^2/4}f(x_1,\cdots,x_i+\sqrt{\lambda}\zeta,\cdots,x_n)d\zeta\\
   &=&
   \frac{1}{2\sqrt{\pi\lambda}}\int^\infty_{-\infty}f(x_1,\cdots,\eta_i,\cdots,x_n)\exp\left(-\frac{(\eta_i-x_i)^2}{4\lambda}\right)d\eta_i.
\end{eqnarray*}
Similarly, (\ref{e2}) is obtained as follows:
$$\exp\left(\lambda\frac{\partial^2}{\partial{x_1^2}}+\cdots+\lambda\frac{\partial^2}{\partial{x_n^2}}\right)f(x)=
\frac{1}{2^n(\lambda\pi)^{n/2}}\int_{\mathbb{R}^n}f(\eta)\exp\left(-\frac{(\eta_1-x_1)^2+\cdots+(\eta_n-x_n)^2}{4\lambda}\right)d\eta.$$
So Example 2.2 is proved.

\textbf{Example 2.3.} Let $f(x),x\in\mathbb{R}^n,\,\lambda>0$, then we have
\begin{equation}\label{yb0}
f(\partial_x)e^{\lambda|x|^2}=
e^{\lambda|x|^2}\frac{1}{2^n(\lambda\pi)^{n/2}}\int_{\mathbb{R}^n}f(\eta)\exp\left(-\frac{|\eta-2\lambda{x}|^2}{4\lambda}\right)d\eta.
\end{equation}

\textbf{Proof.} In Corollary 2.2, let
$g(x)=e^{\lambda|x|^2}\in{C^\infty}(\mathbb{R}^n),\;\lambda>0$, using (\ref{e0}), (\ref{e1}) and (\ref{e2}) we have
\begin{eqnarray*}
% \nonumber to remove numbering (before each equation)
  f(\partial_x)e^{\lambda|x|^2}&=& \left.\exp(\lambda|\partial_\xi|^2)(e^{\xi{x}}f(\xi))\right|_{\xi=0}=
  \left.\exp(\lambda|x+\partial_\xi|^2)f(\xi)\right|_{\xi=0}\\
   &=& e^{\lambda|x|^2}\left.\exp(2\lambda{x}\partial_\xi)\exp(\lambda|\partial_\xi|^2)f(\xi)\right|_{\xi=0}\\
    &=& e^{\lambda|x|^2}\frac{1}{2^n(\lambda\pi)^{n/2}}
 \left.\exp(2\lambda{x}\partial_\xi)\int_{\mathbb{R}^n}f(\eta)\exp\left(-\frac{|\eta-\xi|^2}{4\lambda}\right)d\eta\right|_{\xi=0}\\
    &=&
e^{\lambda|x|^2}\frac{1}{2^n(\lambda\pi)^{n/2}}\left.\int_{\mathbb{R}^n}f(\eta)\exp\left(-\frac{|\eta-\xi-2\lambda{x}|^2}{4\lambda}\right)d\eta\right|_{\xi=0},
\end{eqnarray*}
where
$x\partial_\xi=\langle{x,\partial_\xi}\rangle=x_1\partial_{\xi_1}+x_2\partial_{\xi_2}+\cdots+x_n\partial_{\xi_n}$,
$|x|^2=x_1^2+x_2^2+\cdots+x_n^n$. So we get (\ref{yb0}). Example 2.3 is proved.

Therefore, Corollary 2.2 is an important way to extend the symbols of abstract operators. In fact, limitations on the symbols
$f(t,\xi)$ and the operator element $\partial_x$ will be broadened, and more
general abstract operators will be introduced to deal with problems
involved in our further discussion.

\subsection{Abstract operators and partial differential operators}
\noindent

In order to achieve self consistency theoretically, firstly let us derive the
rules of differentiation from the view of abstract operators, making
the concept of ordinary or partial differential operators on the
basis of the abstract operators itself.

\textbf{Corollary 2.3.} The abstract operators is a kind of generalized partial differential operators. In other words, the partial differential operators $\partial_{x_i}$ and $\partial_x^\alpha$ can be defined as the following abstract operators, namely
\[f_i\left(x_i\frac{\partial}{\partial{x_i}}\right)x^\alpha=\alpha_ix^\alpha\quad\mbox{and}\quad
f(\partial_x)e^{\xi{x}}=\xi^\alpha{e}^{\xi{x}}.\]
Where their symbols $f_i(\alpha_i)=\alpha_i,\,f(\xi)=\xi^\alpha,\,\alpha\in\mathbb{N}^n,\;\xi\in\mathbb{R}_n,\,x\in\mathbb{R}^n$.

\textbf{Proof.} The major parts of rules of differentiation are the
derivative principle of function product, the derivative principle
of compound function and the chain rule of multivariate function.

Let $a,b\in\mathbb{R}$ be the characters of base functions, as the
characteristic  equation, the binomial formula of integer power can
be expressed as
\[{b^n}=\sum^n_{j=0}(-1)^j{n\choose{j}}a^j(a+b)^{n-j},\quad\forall{n}\in\mathbb{N}.\]
By base functions $e^{ax},e^{bx},\,x\in\mathbb{R}^1$ we can combine the binomial
formula with the following operator equality:
\[e^{ax}\frac{d^n}{dx^n}e^{bx}=\sum^n_{j=0}(-1)^j{n\choose{j}}\frac{d^{n-j}}{dx^{n-j}}\left(e^{bx}\frac{d^j}{dx^j}e^{ax}\right).\]
According to the Analytic continuous fundamental theorem, we get
\begin{equation}\label{1*}
v\frac{d^nu}{dx^n}=\sum^n_{j=0}(-1)^j{n\choose{j}}\frac{d^{n-j}}{dx^{n-j}}\left(u\frac{d^jv}{dx^j}\right),\quad\forall{v,u}\in{C^n}(\Omega).
\end{equation}

Similarly the Leibniz rule is tenable, namely
\begin{equation}\label{1}
\frac{d^n}{dx^n}(vu)=\sum^n_{j=0}{n\choose{j}}\frac{d^jv}{dx^j}\frac{d^{n-j}u}{dx^{n-j}},\quad\forall{v,u}\in{C^n}(\Omega).
\end{equation}

Let $n=1$ in (\ref{1}), then for $v=f(x)$
and $u=g(x)$ we have
\begin{equation}\label{3}
    \frac{d}{dx}(f(x)g(x))=f(x)\frac{d}{dx}g(x)+g(x)\frac{d}{dx}f(x),\quad\forall{f,g}\in{C}^1(\Omega).
\end{equation}

Let $\varphi(x)=f_1(x)f_2(x)\cdots{f}_n(x)$, generally we have
\[\frac{d}{dx}\varphi(x)=\sum^n_{j=1}\prod^n_{i=1\atop
i\neq{j}}f_i(x)\frac{d}{dx}f_j(x),\quad\forall{f_j(x)}\in{C}^1(\Omega),\;j=1,2,\cdots,n.\]

If $f_1(x)=f_2(x)=\cdots=f_n(x)=y=g(x)\in{C}^1(\Omega)$, then for differentiable
function $y=g(x)$ we have
\[\frac{d}{dx}y^n=ny^{n-1}\frac{dy}{dx}\quad\mbox{or}\quad\frac{d}{dx}y^n=\frac{dy}{dx}\frac{d}{dy}y^n,\qquad\forall{n}\in\mathbb{N}^1.\]
According to the Analytic continuous fundamental theorem, we have
the derivative principle of compound function
$f(y)\in{C^1}(V),\,y\in{V}\subset\mathbb{R}^1$:
\begin{equation}\label{4}
\frac{d}{dx}f(g(x))=\frac{dy}{dx}\frac{d}{dy}f(y),\qquad{y}=g(x)\in{C}^1(\Omega),\quad{x}\in\Omega\subset\mathbb{R}^1.
\end{equation}

By using (\ref{3}) and (\ref{4}), for differentiable functions
$x_1(t)$ and $x_2(t)$, $t\in{I}\subset\mathbb{R}^1$ we have
\[\frac{d}{dt}(x^{\alpha_1}_1x^{\alpha_2}_2)=\frac{dx^{\alpha_1}_1}{dx_1}\frac{dx_1}{dt}x^{\alpha_2}_2
+x^{\alpha_1}_1\frac{dx^{\alpha_2}_2}{dx_2}\frac{dx_2}{dt},\qquad(\alpha_1,\alpha_2)\in\mathbb{N}^2.\]

Taking this one as the characteristic equation, by base function
$x^\alpha$, $x\in\mathbb{R}^2,\;\alpha\in\mathbb{N}^2$, combining it
with the following operator equality:
\[\frac{d}{dt}x^\alpha=\frac{dx_1}{dt}\frac{\partial}{\partial{x}_1}x^\alpha+\frac{dx_2}{dt}\frac{\partial}{\partial{x}_2}x^\alpha,
\qquad{x}(t)\in\mathbb{R}^2\quad{t}\in{I}\subset\mathbb{R}^1\quad\alpha\in\mathbb{N}^2.\]
According to the Analytic continuous fundamental theorem, $\forall{f(x)}\in{C^1}(\Omega),\;x\in\Omega\subset\mathbb{R}^2$ we have
\[\frac{d}{dt}f(x)=\frac{dx_1}{dt}\frac{\partial}{\partial{x}_1}f(x)+\frac{dx_2}{dt}\frac{\partial}{\partial{x}_2}f(x),
\quad x_i(t)\in{C}^1(I),\;i=1,2.\]

Clearly, $\forall{f(x)}\in{C^1}(\Omega),\;x\in\Omega\subset\mathbb{R}^n$, $t\in{I}\subset\mathbb{R}^1$ we can
generally derive
\begin{equation}\label{5}
\frac{d}{dt}f(x)=\frac{dx_1}{dt}\frac{\partial}{\partial{x}_1}f(x)+\frac{dx_2}{dt}\frac{\partial}{\partial{x}_2}f(x)
+\cdots+\frac{dx_n}{dt}\frac{\partial}{\partial{x}_n}f(x),\quad x_i(t)\in{C}^1(I),\;i=1,2,\cdots,n.
\end{equation}

If taking $n\in\mathbb{N}$ as the character  of base function, the
binomial formula can be expressed as the following characteristic
equation:
\[(x+h)^n=\sum^\infty_{j=0}\frac{h^j}{j!}\frac{d^j}{dx^j}x^n,\quad{x}\in\mathbb{C}^1.\]
According to the Analytic continuous fundamental theorem, the Taylor formula is
tenable for any analytic function, namely
\begin{equation}\label{2}
f(x+h)=\sum^\infty_{j=0}\frac{h^j}{j!}\frac{d^j}{dx^j}f(x),\quad\forall f(x)\in{C^\infty}(\Omega),\;x\in\Omega\subset\mathbb{C}^1,\;|h|<R.
\end{equation}

Let $a,b\in\mathbb{R}$ and $a\leq{b}$, without losing the universality,
assuming $a\geq0,\;b\geq0$ we have
\[na^{n-1}\leq(a^{n-1}+a^{n-2}b+a^{n-3}b^2+\cdots+ab^{n-2}+b^{n-1})\leq{n}b^{n-1},\quad{n}=1,2,\cdots.\]
namely
\[na^{n-1}\leq\frac{b^n-a^n}{b-a}\leq{n}b^{n-1}\quad\mbox{or}\quad
{a}\leq\sqrt[n-1]{\frac{1}{n}\frac{b^n-a^n}{b-a}}\leq{b}.\]
Therefore, if $a\leq{b}$, then we have $a\leq{c}\leq{b}$, making
\[nc^{n-1}=\frac{b^n-a^n}{b-a}\quad\mbox{or}\quad\frac{d}{dc}c^n=\frac{b^n-a^n}{b-a},\qquad\forall{n}\in\mathbb{N}^1.\]
Taking this one as the characteristic equation, according to the
Analytic continuous fundamental theorem we get the Lagrange mean
value theorem, namely

If $a\leq{b}$, then $\exists{c}\in[a,b]$ making
\begin{equation}\label{lg}
\frac{d}{dc}f(c)=\frac{f(b)-f(a)}{b-a},\quad\forall{f(x)}\in{C^1}[a,b].
\end{equation}
So the Corollary 2.3 is proved.

\textbf{Example 2.4.} Let $y=x^2,\,x\in\mathbb{R}^1$, then $\forall{f(y)}\in{C^k}(\Omega),\;y\in\Omega\subset\mathbb{R}^1$, we have
\begin{equation}\label{y5}
\frac{d^k}{dx^k}f(x^2)=\sum^{[k/2]}_{j=0}\frac{k!}{j!\,(k-2j)!}(2x)^{k-2j}\frac{d^{k-j}}{dy^{k-j}}f(y).
\end{equation}
If $y=x^3,\,x\in\mathbb{R}^1$, then $\forall{f(y)}\in{C^n}(\Omega)$
\begin{equation}\label{y5.3}
\frac{d^n}{dx^n}f(x^3)=\sum^{[n/2]}_{k=0}\sum^k_{j=0}\frac{3^{n-k-2j}n!}{j!\,(k-j)!\,(n-2k-j)!}x^{2n-3k-3j}\frac{d^{n-k-j}}{dy^{n-k-j}}f(y).
\end{equation}

\textbf{Proof.} Let $x\in\mathbb{R}^1,\,\xi\in\mathbb{R}_1,\,\zeta=3x^2\lambda$, $\lambda$ be a real- or complex-valued, applying Corollary 2.2 and Example 2.1 we have
\begin{eqnarray*}
  f(\partial_x)e^{\lambda{x^3}} &=& \left.\exp\left(\lambda\frac{\partial^3}{\partial\xi^3}\right)(e^{\xi{x}}f(\xi))\right|_{\xi=0}
  =\left.\exp\left(\lambda\left(x+\frac{\partial}{\partial\xi}\right)^3\right)f(\xi)\right|_{\xi=0}\\
   &=& \left.\exp\left(\lambda{x^3}+3x^2\lambda\frac{\partial}{\partial\xi}+3x\lambda\frac{\partial^2}{\partial\xi^2}+\lambda\frac{\partial^3}{\partial\xi^3}\right)f(\xi)\right|_{\xi=0}\\
   &=& e^{\lambda{x^3}}\left.\exp\left(3x\lambda\frac{\partial^2}{\partial\xi^2}+\lambda\frac{\partial^3}{\partial\xi^3}\right)f(\xi+3x^2\lambda)\right|_{\xi=0}
   =e^{\lambda{x^3}}\exp\left(3x\lambda\frac{\partial^2}{\partial\zeta^2}+\lambda\frac{\partial^3}{\partial\zeta^3}\right)f(\zeta).
\end{eqnarray*}
So $\forall{f(\zeta)}\in{C^\infty}(\Omega),\;\zeta\in\Omega\subset\mathbb{R}^1\;\mbox{or}\;\mathbb{C}^1$ we get
\begin{equation}\label{y5.4}
e^{-\lambda{x^3}}f(\partial_x)e^{\lambda{x^3}}=\exp\left(3x\lambda\frac{\partial^2}{\partial\zeta^2}+\lambda\frac{\partial^3}{\partial\zeta^3}\right)f(\zeta),
\end{equation}
where $x\in\mathbb{R}^1,\;\zeta=3x^2\lambda$. We can write (\ref{y5.4}) as
\begin{equation}\label{y5.5}
f\left(\frac{d}{dx}\right)e^{\lambda{x^3}}=\sum^\infty_{k=0}\frac{\lambda^k}{k!}\left(3x+\frac{\partial}{\partial\zeta}\right)^kf^{(2k)}(\zeta)\,e^{\lambda{x^3}},\quad|\lambda|<R.
\end{equation}
Let $y=x^3,\,f(\zeta)=\zeta^n,\,n\in\mathbb{N}$, then (\ref{y5.5}) gives following operator equality:
\[\frac{d^n}{dx^n}\,e^{\lambda{x^3}}=\sum^{[n/2]}_{k=0}\sum^k_{j=0}\frac{3^{n-k-2j}n!}{j!\,(k-j)!\,(n-2k-j)!}x^{2n-3k-3j}\frac{d^{n-k-j}}{dy^{n-k-j}}\,e^{\lambda{y}}.\]
So we get (\ref{y5.3}) according to the Analytic continuous fundamental theorem. Similarly we get (\ref{y5}).

Therefore, the mathematical concept of an abstract operators, generalizes the notion of partial differential operators.

\subsection{Algorithms of abstract operators}
\noindent

According to the Analytic continuous fundamental theorem and
Definition 2.4, by using the algebraic properties of two symbols $\cos(h\xi)$ and $\sin(h\xi)$,
Guang-Qing Bi has obtained a series of algorithms of abstract
operators similar to (\ref{3})-(\ref{5}):

\textbf{Theorem BI1.} (See \cite{bi97}) Let $x\in\mathbb{R}^n,\;h\in\mathbb{R}_n$,
$h\partial_x=\langle{h,\partial_x}\rangle$, for the abstract operators $\cos(h\partial_x)\,\mbox{and}\,\sin(h\partial_x)$ we have
\begin{equation}\label{y0}
  \cos(h\partial_x)f(x)=\Re[f(x+ih)],\quad\sin(h\partial_x)f(x)=\Im[f(x+ih)],
\end{equation}
$\forall{f(z)}\in{C}^\infty(\Omega),\,z=x+iy\in\Omega\subset\mathbb{C}^n$.

\parbox{13cm}{\begin{eqnarray*}\label{y1}
                \sin(h\partial_x)(uv) &=&
                \cos(h\partial_x)v\cdot\sin(h\partial_x)u+\sin(h\partial_x)v\cdot\cos(h\partial_x)u,\\
                \cos(h\partial_x)(uv) &=&
                \cos(h\partial_x)v\cdot\cos(h\partial_x)u-\sin(h\partial_x)v\cdot\sin(h\partial_x)u.
              \end{eqnarray*}}\hfill\parbox{1cm}{\begin{eqnarray}\end{eqnarray}}

\parbox{13cm}{\begin{eqnarray*}\label{y2}
                \sin(h\partial_x)\frac{u}{v} &=& \frac{\cos(h\partial_x)v\cdot\sin(h\partial_x)u-\sin(h\partial_x)v\cdot\cos(h\partial_x)u}
                {(\cos(h\partial_x)v)^2+(\sin(h\partial_x)v)^2}, \\
                \cos(h\partial_x)\frac{u}{v} &=& \frac{\cos(h\partial_x)v\cdot\cos(h\partial_x)u+\sin(h\partial_x)v\cdot\sin(h\partial_x)u}
                {(\cos(h\partial_x)v)^2+(\sin(h\partial_x)v)^2}.
              \end{eqnarray*}}\hfill\parbox{1cm}{\begin{eqnarray}\end{eqnarray}}

\textbf{Proof.} See Reference \cite{bi97} (G.Q. Bi, 1997).

\textbf{Theorem BI2.} (See \cite{bi97}) Let $h_0\in\mathbb{C},\;x(t)\in\mathbb{R}^n,\;t\in\mathbb{R}^1,\;X\in\mathbb{R}^n,\,Y\in\mathbb{R}_n$,
$Y\partial_X=\langle{Y,\partial_X}\rangle=Y_1\partial_{X_1}+\cdots+Y_n\partial_{X_n}$,
then we have

\parbox{13cm}{\begin{eqnarray*}\label{y3}
 \sin(h_0\partial_t)f(x(t)) &=& \sin(Y\partial_X)f(X),\\
 \cos(h_0\partial_t)f(x(t)) &=& \cos(Y\partial_X)f(X).
\end{eqnarray*}}\hfill\parbox{1cm}{\begin{eqnarray}\end{eqnarray}}
Where $X_j=\cos(h_0\partial_t)x_j(t),\;Y_j=\sin(h_0\partial_t)x_j(t),\;j=1,\cdots,n$.

For example, when $n=1$, the (\ref{y3}) can easily be restated as

\parbox{13cm}{\begin{eqnarray*}\label{yb1}
\sin\left(h_0\frac{d}{dt}\right)f(x(t)) &=& \sin\left(Y\frac{\partial}{\partial{X}}\right)f(X),\\
\cos\left(h_0\frac{d}{dt}\right)f(x(t)) &=& \cos\left(Y\frac{\partial}{\partial{X}}\right)f(X).
\end{eqnarray*}}\hfill\parbox{1cm}{\begin{eqnarray}\end{eqnarray}}
Where
$Y=\sin(h_0\partial_t)x(t),\;X=\cos(h_0\partial_t)x(t),\;t\in\mathbb{R}^1,\;h_0\in\mathbb{C}$.

If $n=2$, then (\ref{y3}) can easily be restated as

\parbox{13cm}{\begin{eqnarray*}\label{yb2}
\sin\left(h_0\frac{d}{dt}\right)f(x(t),y(t)) &=& \sin\left(Y_x\frac{\partial}{\partial{X}_x}+Y_y\frac{\partial}{\partial{X}_y}\right)f(X_x,X_y),\\
\cos\left(h_0\frac{d}{dt}\right)f(x(t),y(t)) &=& \cos\left(Y_x\frac{\partial}{\partial{X}_x}+Y_y\frac{\partial}{\partial{X}_y}\right)f(X_x,X_y).
\end{eqnarray*}}\hfill\parbox{1cm}{\begin{eqnarray}\end{eqnarray}}
Where
$Y_x=\sin(h_0\partial_t)x(t),\;X_x=\cos(h_0\partial_t)x(t),\;Y_y=\sin(h_0\partial_t)y(t),\;X_y=\cos(h_0\partial_t)y(t)$.

\textbf{Theorem 2.1.}  Let $u=g(y)$ be a monotonic function on its
domain, if $y=f(bx)$ is the inverse function of $bx=g(y)$, namely
$g(f(bx))=bx$, where
$bx=b_1x_1+b_2x_2+\cdots+b_nx_n,\;bh=b_1h_1+b_2h_2+\cdots+b_nh_n$,
then $\sin(h\partial_x)f(bx)$ (denoted by $Y$) and
$\cos(h\partial_x)f(bx)$ (denoted by $X$) can be determined by the
following set of equations:
\begin{equation}\label{y4}
    \left\{\begin{array}{l@{\qquad}l}\displaystyle
    \cos\left(Y\frac{\partial}{\partial{X}}\right)g(X)=bx,&x\in\mathbb{R}^n,\;b\in\mathbb{R}_n,\\\displaystyle
    \sin\left(Y\frac{\partial}{\partial{X}}\right)g(X)=bh, &h\in\mathbb{R}^n.
    \end{array}\right.
\end{equation}

\textbf{Proof.} According to the Analytic continuous fundamental
theorem and Definition 2.4, we have
\[e^{ih\partial_x}g(y)=g(f(bx+ibh))=g(e^{ih\partial_x}f(bx))=g(X+iY),\]
and considering operator algebras of abstract operators we obtain

\parbox{13cm}{\begin{eqnarray*}\label{ey3}
                \sin(h\partial_x)g(y) &=& \sin\left(Y\frac{\partial}{\partial{X}}\right)g(X),\\
                \cos(h\partial_x)g(y) &=& \cos\left(Y\frac{\partial}{\partial{X}}\right)g(X).
              \end{eqnarray*}}
              \hfill\parbox{1cm}{\begin{eqnarray}\end{eqnarray}}
Where
$y=f(bx),\,bx=\langle{b,x}\rangle,\;X=\cos(h\partial_x)f(bx),\;Y=\sin(h\partial_x)f(bx),\;x\in\mathbb{R}^n,\;b\in\mathbb{R}_n,\;h\in\mathbb{R}^n$.

In the left side of (\ref{ey3}), when $u=g(y)$ is a monotonic function on its
domain, and there is function $y=f(bx)$ making $g(f(bx))=bx$, then we have respectively
$\sin(h\partial_x)g(y)=\sin(h\partial_x)bx=bh$ and $\cos(h\partial_x)g(y)=\cos(h\partial_x)bx=bx$.
So Theorem 2.1 is proved.

\textbf{Example 2.5.} Let $f(x)\in{L}^2([-c,c])$ be the square wave function on $\mathbb{R}^1$ defined by
\[f(x)=\left\{\begin{array}{r@{\qquad}l}
+1,&2Kc-c/2<x<2Kc+c/2\\
-1,&2Kc+c/2<x<2Kc+3c/2,\end{array}\right.\]
where $K=0,\pm1,\pm2,\cdots$, then $f(x)$ can be expressed as
\[f(x)=\frac{4}{\pi}\left.\cos\left(\frac{\pi{x}}{c}\frac{\partial}{\partial{z}}\right)\arctan{e^z}\right|_{z=0}.\]

\textbf{Proof.} (See \cite{bigy}) For $\sin(h\partial_x)\arctan{bx}$ and
$\cos(h\partial_x)\arctan{bx}$, using (\ref{y0}), (\ref{y2}) and (\ref{y4}) we can obtain

\parbox{13cm}{\begin{eqnarray*}\label{ly8}
                \sin(h\partial_x)\arctan{bx} &=& \frac{1}{2}\textrm{tanh}^{-1}\frac{2bh}{1+(bx)^2+(bh)^2}, \\
                \cos(h\partial_x)\arctan{bx} &=&
                \frac{1}{2}\arctan\frac{2bx}{1-(bx)^2-(bh)^2},
              \end{eqnarray*}}\hfill\parbox{1cm}{\begin{eqnarray}\end{eqnarray}}
where
$bx=b_1x_1+b_2x_2+\cdots+b_nx_n,\;bh=b_1h_1+b_2h_2+\cdots+b_nh_n$.

By (\ref{yb1}) and (\ref{ly8}) we have
\begin{eqnarray*}
% \nonumber to remove numbering (before each equation)
   f(x)&=& \frac{4}{\pi}\left.\cos\left(\frac{\pi{x}}{c}\frac{\partial}{\partial{z}}\right)\arctan{e^z}\right|_{z=0}=\,
   \frac{4}{\pi}\left.\cos\left(Y\frac{\partial}{\partial{X}}\right)\arctan{X}\right|_{z=0}\\
   &=&\left.\frac{2}{\pi}\arctan\frac{2X}{1-(X^2+Y^2)}\right|_{z=0}\,=\,
   \frac{2}{\pi}\arctan\frac{2\cos(\pi{x}/c)}{1-\left(\cos^2(\pi{x}/c)+\sin^2(\pi{x}/c)\right)}.
\end{eqnarray*}
When $\cos(\pi{x}/c)\neq0$ or $|x|\neq{kc}+c/2,\;k=0,1,2,\cdots$,
the above expression turns into
\[\frac{4}{\pi}\left.\cos\left(\frac{\pi{x}}{c}\frac{\partial}{\partial{z}}\right)\arctan{e^z}\right|_{z=0}=\left\{\begin{array}{r@{\qquad}l}
(2/\pi)\arctan(+\infty)=+1,&2Kc-c/2<x<2Kc+c/2\\
(2/\pi)\arctan(-\infty)=-1,&2Kc+c/2<x<2Kc+3c/2,\end{array}\right.\]
where $K=0,\pm1,\pm2,\cdots$.

More generally, we have

\textbf{Theorem 2.2.} Let ${f(x)}\in{L}^2([-l,l]),\;f(x+2l)=f(x)$, if
$S_+(t)$ and $S_-(t)$ are respectively given by
\begin{equation}\label{f1'}
S_+(t)=\frac{1}{2l}\int^l_{-l}f(\xi)\frac{1-t^2}{1-2t\cos(\pi\xi/l)+t^2}d\xi,
\end{equation}
\begin{equation}\label{f2'}
 S_-(t)=\frac{1}{l}\int^l_{-l}f(\xi)\frac{t\sin(\pi\xi/l)}{1-2t\cos(\pi\xi/l)+t^2}d\xi,
\end{equation}
then $\exists{f_z(x)}\in{C^\infty}(\mathbb{R}^1)$ with the following
form
\begin{equation}\label{f3'}
f_z(x)=\cos\left(\frac{\pi{x}}{l}\frac{\partial}{\partial{z}}\right)S_+(e^z)+\sin\left(\frac{\pi{x}}{l}\frac{\partial}{\partial{z}}\right)S_-(e^z),\quad
-\infty<z<0
\end{equation}
making
$$\left.f_z(x)\right|_{z=0}=\lim_{z\rightarrow0^-}f_z(x)\rightharpoonup{f(x)},\quad\forall{f(x)}\in{L}^2([-l,l]).$$

\textbf{Proof.} By the algorithms (\ref{y2}) we have
\[\cos\left(\frac{\pi\xi}{l}\frac{\partial}{\partial{z}}\right)\frac{1+e^z}{1-e^z}=\frac{1-e^{2z}}{1-2e^{z}\cos(\pi\xi/l)+e^{2z}}.\]
\[\sin\left(\frac{\pi\xi}{l}\frac{\partial}{\partial{z}}\right)\frac{e^z}{1-e^z}=\frac{e^z\sin(\pi\xi/l)}{1-2e^{z}\cos(\pi\xi/l)+e^{2z}}.\]
So $f_z(x)$ can be expressed as
\begin{eqnarray*}
% \nonumber to remove numbering (before each equation)
   f_z(x) &=&
   \cos\left(\frac{\pi{x}}{l}\frac{\partial}{\partial{z}}\right)\left[\frac{1}{2l}\int^l_{-l}f(\xi)
   \cos\left(\frac{\pi\xi}{l}\frac{\partial}{\partial{z}}\right)\frac{1+e^z}{1-e^z}d\xi\right]\\
   & & +\,\sin\left(\frac{\pi{x}}{l}\frac{\partial}{\partial{z}}\right)\left[\frac{1}{l}\int^l_{-l}f(\xi)
   \sin\left(\frac{\pi\xi}{l}\frac{\partial}{\partial{z}}\right)\frac{e^z}{1-e^z}d\xi\right] \\
   &=& \frac{1}{2l}\int^l_{-l}f(\xi)\cos\left(\frac{\pi{x}}{l}\frac{\partial}{\partial{z}}\right)
   \cos\left(\frac{\pi\xi}{l}\frac{\partial}{\partial{z}}\right)\left[1+\frac{2e^z}{1-e^z}\right]d\xi\\
   & & +\,\frac{1}{l}\int^l_{-l}f(\xi)\sin\left(\frac{\pi{x}}{l}\frac{\partial}{\partial{z}}\right)
   \sin\left(\frac{\pi\xi}{l}\frac{\partial}{\partial{z}}\right)\frac{e^z}{1-e^z}d\xi = \frac{1}{2l}\int^l_{-l}f(\xi)d\xi\\
   & & +\,\frac{1}{l}\int^l_{-l}f(\xi)\left[\cos\left(\frac{\pi{x}}{l}\frac{\partial}{\partial{z}}\right)
   \cos\left(\frac{\pi\xi}{l}\frac{\partial}{\partial{z}}\right)+\sin\left(\frac{\pi{x}}{l}\frac{\partial}{\partial{z}}\right)
   \sin\left(\frac{\pi\xi}{l}\frac{\partial}{\partial{z}}\right)\right]\frac{e^z}{1-e^z}d\xi\\
   &=&
   \frac{1}{2l}\int^l_{-l}f(\xi)d\xi+\frac{1}{l}\int^l_{-l}f(\xi)\cos\left(\frac{\pi(x-\xi)}{l}\frac{\partial}{\partial{z}}\right)\frac{e^z}{1-e^z}d\xi\\
   &=& \frac{1}{2l}\int^l_{-l}f(\xi)\cos\left(\frac{\pi(x-\xi)}{l}\frac{\partial}{\partial{z}}\right)\frac{1+e^z}{1-e^z}d\xi\\
   &=& \frac{1}{2l}\int^l_{-l}f(\xi)\frac{1-e^{2z}}{1-2e^{z}\cos(\pi(x-\xi)/l)+e^{2z}}d\xi,\quad-\infty<z<0.
\end{eqnarray*}
\begin{eqnarray*}
 \lim_{z\rightarrow0^-}f_z(x) &=& \lim_{z\rightarrow0^-}\frac{1}{2l}\int^l_{-l}f(\xi)\frac{1-e^{2z}}{1-2e^{z}\cos(\pi(x-\xi)/l)+e^{2z}}d\xi\\
 &=&
\int^l_{-l}f(\xi)\lim_{z\rightarrow0^-}\frac{1}{2l}\frac{1-e^{2z}}{1-2e^{z}\cos(\pi(x-\xi)/l)+e^{2z}}d\xi
\rightharpoonup\int^l_{-l}f(\xi)\delta(x-\xi)d\xi =f(x).
\end{eqnarray*}
Thus Theorem 2.2 is proved.

In Theorem 2.2, $f(x)=\lim_{z\rightarrow0^-}f_z(x)\in{L}^2([-l,l])$, which may be discontinuous on real axis, such as Example 2.5, $f(x)$ is the square wave function, however, $f_z(x)$ corresponding to $f(x)$ is a $C^\infty$ function. Therefore, puppose that $g(\partial_x)$ is an abstract operators, if $\lim_{z\rightarrow0^-}g(\partial_x)f_z(x)\in{L}^2([-l,l])$, then $g(\partial_x)f(x)$ makes sense in a broad sense, which can be defined as
\begin{equation}\label{f4'}
    g(\partial_x)f(x)=\lim_{z\rightarrow0^-}g(\partial_x)f_z(x).
\end{equation}

\textbf{Example 2.6.} (See \cite{bigy})
Let
$S(t)=\sum^\infty_{n=0}a_nt^n,\;t\in\mathbb{R}^1,\;0\leq{t}\leq{r},\;0<r<+\infty$,
if
\[\sum^\infty_{n=0}a_n\cos\frac{n\pi{x}}{c}=\left.\cos\left(\frac{\pi{x}}{c}\frac{\partial}{\partial{z}}\right)S(e^z)\right|_{z=0},\]
\[\sum^\infty_{n=0}a_n\sin\frac{n\pi{x}}{c}=\left.\sin\left(\frac{\pi{x}}{c}\frac{\partial}{\partial{z}}\right)S(e^z)\right|_{z=0}\]
are tenable on $a<x<b,\;x\in\mathbb{R}^1$, then
$\forall{x_0}\in[0,\frac{b-a}{2})$,

\parbox{13cm}{\begin{eqnarray*}\label{yb5}
 \sum^\infty_{n=0}a_n\cos\frac{n\pi{x_0}}{c}\cos\frac{n\pi{x}}{c}&=&
 \cosh\left(x_0\frac{\partial}{\partial{x}}\right)\left.\cos\left(\frac{\pi{x}}{c}\frac{\partial}{\partial{z}}\right)S(e^z)\right|_{z=0},\\
 \sum^\infty_{n=0}a_n\cos\frac{n\pi{x_0}}{c}\sin\frac{n\pi{x}}{c}&=&
 \cosh\left(x_0\frac{\partial}{\partial{x}}\right)\left.\sin\left(\frac{\pi{x}}{c}\frac{\partial}{\partial{z}}\right)S(e^z)\right|_{z=0}
\end{eqnarray*}}\hfill\parbox{1cm}{\begin{eqnarray}\end{eqnarray}}
are tenable on $a+x_0<x<b-x_0$. Accordingly, for any definite value
$x\in(a,b)$, $x_0$ in (\ref{yb5}) takes values in the following interval
\[0<x_0<\left\{\begin{array}{lr}x-a,&a<x\leq(a+b)/2,\\b-x,&b>x\geq(a+b)/2.\end{array}\right.\]
Where
\[\cos\frac{n\pi{x_0}}{c}\cos\frac{n\pi{x}}{c}=\cosh\left(x_0\frac{\partial}{\partial{x}}\right)\cos\frac{n\pi{x}}{c},\]
\[\cos\frac{n\pi{x_0}}{c}\sin\frac{n\pi{x}}{c}=\cosh\left(x_0\frac{\partial}{\partial{x}}\right)\sin\frac{n\pi{x}}{c}.\]
Therefore, $\cosh(x_0\partial_x)$ and $\sinh(x_0\partial_x)$ are the continuous operators on $L^2([-c,c])$.

\textbf{Theorem 2.3.}
 Let $t>0,\,\Re(s)>0,\;\forall{f(x)}\in\mathscr{S}(\mathbb{R}^1)$, if $F_+(s)$ and $F_-(s)$ are
 respectively given by
\begin{equation}\label{f1}
F_+(s)=\frac{1}{\pi}\mathscr{L}\int^\infty_{-\infty}f(\xi)\cos(t\xi)d\xi,\qquad
F_-(s)=\frac{1}{\pi}\mathscr{L}\int^\infty_{-\infty}f(\xi)\sin(t\xi)d\xi,
\end{equation}
or
\begin{equation}\label{f1s}
F_+(s)=\frac{1}{\pi}\int^\infty_{-\infty}f(\xi)\frac{s}{\xi^2+s^2}d\xi,\qquad
F_-(s)=\frac{1}{\pi}\int^\infty_{-\infty}f(\xi)\frac{\xi}{\xi^2+s^2}d\xi,
\end{equation}
where $\mathscr{L}$ is the Laplace transform operator, then we have
\begin{equation}\label{f3}
f(x)=\left.\cos\left(x\frac{\partial}{\partial{s}}\right)F_+(s)\right|_{s=0}+\left.\sin\left(-x\frac{\partial}{\partial{s}}\right)F_-(s)\right|_{s=0}.
\end{equation}

\textbf{Proof.} By (\ref{t4}) we have
\[\mathscr{L}\cos(t\xi)=\cos\left(\xi\frac{\partial}{\partial{s}}\right)\frac{1}{s},\quad
\mathscr{L}\sin(t\xi)=\sin\left(-\xi\frac{\partial}{\partial{s}}\right)\frac{1}{s},\;\Re(s)>0.\]
So by (\ref{f1}) and (\ref{f3}),
$\forall{f(x)}\in\mathscr{S}(\mathbb{R}^1),\,\exists{f_s(x)}\in{C^\infty}(\mathbb{R}^1)$
can be expressed as
\begin{eqnarray*}
% \nonumber to remove numbering (before each equation)
   f_s(x) &=& \cos\left(x\frac{\partial}{\partial{s}}\right)F_+(s)+\sin\left(-x\frac{\partial}{\partial{s}}\right)F_-(s) \\
   &=& \cos\left(x\frac{\partial}{\partial{s}}\right)\left[\frac{1}{\pi}\int^\infty_{-\infty}f(\xi)\mathscr{L}\cos(t\xi)d\xi\right]
   +\sin\left(-x\frac{\partial}{\partial{s}}\right)\left[\frac{1}{\pi}\int^\infty_{-\infty}f(\xi)\mathscr{L}\sin(t\xi)d\xi\right]\\
   &=& \frac{1}{\pi}\cos\left(x\frac{\partial}{\partial{s}}\right)\int^\infty_{-\infty}f(\xi)
   \cos\left(\xi\frac{\partial}{\partial{s}}\right)\frac{1}{s}\,d\xi
   +\frac{1}{\pi}\sin\left(-x\frac{\partial}{\partial{s}}\right)\int^\infty_{-\infty}f(\xi)
   \sin\left(-\xi\frac{\partial}{\partial{s}}\right)\frac{1}{s}\,d\xi\\
   &=& \frac{1}{\pi}\int^\infty_{-\infty}f(\xi)\left[\cos\left(x\frac{\partial}{\partial{s}}\right)
   \cos\left(\xi\frac{\partial}{\partial{s}}\right)+\sin\left(-x\frac{\partial}{\partial{s}}\right)
   \sin\left(-\xi\frac{\partial}{\partial{s}}\right)\right]\frac{1}{s}\,d\xi\\
   &=& \frac{1}{\pi}\int^\infty_{-\infty}f(\xi)\cos\left((x-\xi)\frac{\partial}{\partial{s}}\right)\frac{1}{s}\,d\xi
   =\frac{1}{\pi}\int^\infty_{-\infty}f(\xi)\frac{s}{(x-\xi)^2+s^2}d\xi,\;\Re(s)>0.
\end{eqnarray*}
Therefore, we get
\begin{eqnarray*}
 \lim_{s\rightarrow0^+}f_s(x) &=& \lim_{s\rightarrow0^+}\frac{1}{\pi}\int^\infty_{-\infty}f(\xi)\frac{s}{(x-\xi)^2+s^2}d\xi\\
 &=&
\int^\infty_{-\infty}f(\xi)\lim_{s\rightarrow0^+}\frac{1}{\pi}\frac{s}{(x-\xi)^2+s^2}d\xi\rightharpoonup\int^\infty_{-\infty}f(\xi)\delta(x-\xi)d\xi
=f(x).
\end{eqnarray*}
Thus Theorem 2.3 is proved.

Substituting (\ref{f3}) into (\ref{f1}), then Theorem 2.3 gives

\textbf{Corollary 2.4.}  Let $t>0,\,F(s)=\mathscr{L}f(t)$, then $\forall{f(t)}\in\mathscr{S}(\mathbb{R}^1)$ we have
\begin{equation}\label{f1+}
\mathscr{L}^{-1}F(s)=\frac{2}{\pi}\int^\infty_0 \!\left.\cos\left(\xi\frac{\partial}{\partial{s}}\right)F(s)\right|_{s=0}\cos(t\xi)d\xi,
\end{equation}
and
\begin{equation}\label{f1++}
\mathscr{L}^{-1}F(s)=\frac{2}{\pi}\int^\infty_0\!\left.\sin\left(-\xi\frac{\partial}{\partial{s}}\right)F(s)\right|_{s=0}\sin(t\xi)d\xi.
\end{equation}

In fact, if $f(x)$ is the rational proper franctions, then the integral (\ref{f1}) are absolutely comvergent. Therefore, by (\ref{f1}) and (\ref{f3}) we get the following Corollary:

\textbf{Corollary 2.5.} Let $f(x),\,x\in\mathbb{R}^1$ be the rational proper franctions, if analytic functions $F_+(s)$ satisfies the following operator equation
\begin{equation}\label{f3+}
\left.\cos\left(x\frac{\partial}{\partial{s}}\right)F_+(s)\right|_{s=0}=\frac{1}{2}(f(x)+f(-x)),
\end{equation}
then we have
\begin{equation}\label{f1+3}
\int^\infty_{-\infty}f(x)\cos(tx)dx=\pi\mathscr{L}^{-1}F_+(s),\quad t>0.
\end{equation}
If analytic functions $F_-(s)$ satisfies the following operator equation
\begin{equation}\label{f3++}
\left.\sin\left(-x\frac{\partial}{\partial{s}}\right)F_-(s)\right|_{s=0}=\frac{1}{2}(f(x)-f(-x)),
\end{equation}
then we have
\begin{equation}\label{f1+3+}
\int^\infty_{-\infty}f(x)\sin(tx)dx=\pi\mathscr{L}^{-1}F_-(s),\quad t>0.
\end{equation}

Obviously, $F_+(s)$ and $F_-(s)$ are the rational proper franctions when $f(x)$ is a rational proper franctions. So we can use the method of undetermined coefficients to solve these two operator equations (\ref{f3+}) and (\ref{f3++}) for a given rational proper franctions $f(x)$.

Because of the introduction of Definition 2.3, it is necessary to add the following algorithms:

\textbf{Theorem 2.4.} Let $a=(a_1,\cdots,a_n)\in\mathbb{C}^n,\,x\in\mathbb{R}^n,\,a^{x\partial_x}=a_1^{x_1\partial_1}a_2^{x_2\partial_2}\cdots{a}_n^{x_n\partial_n}$, then
we have
\begin{equation}\label{15}
    a^{x\partial_x}f(x)=f(a_1x_1,\cdots,a_ix_i,\cdots,a_nx_n),\quad\forall{f}(x)\in{C}^\infty(\Omega).
\end{equation}

\textbf{Proof.} By Definition 2.3 for the base function $x^\alpha,\,x\in\mathbb{R}^n,\,\alpha\in\mathbb{N}^n$ we have
\[a^{x\partial_x}x^\alpha=(a_1x_1)^{\alpha_1}(a_2x_2)^{\alpha_2}\cdots(a_nx_n)^{\alpha_n}.\]
So Theorem 2.4 is proved by the Analytic continuous fundamental theorem.

\textbf{Theorem 2.5.} Let
$X=(\rho_1\cos\theta_1,\cdots,\rho_n\cos\theta_n),\;
Y=(\rho_1\sin\theta_1,\cdots,\rho_n\sin\theta_n)$, introducing the
notation
\[\theta\,\rho\frac{\partial}{\partial\rho}=\langle\theta,\rho\partial_\rho\rangle=\theta_1\,\rho_1\frac{\partial}{\partial\rho_1}+\cdots+
\theta_n\,\rho_n\frac{\partial}{\partial\rho_n},\;\;Y\partial_X=\langle{Y,\partial_X}\rangle=Y_1\frac{\partial}{\partial{X}_1}+\cdots+Y_n\frac{\partial}{\partial{X}_n},\]
then
$\forall{f}(\rho)\in{C}^\infty(\Omega),\;\rho\in\mathbb{R}^n,\;\theta\in\mathbb{R}_n$,
we have

\parbox{11cm}{\begin{eqnarray*}\label{15*3}
             % \nonumber to remove numbering (before each equation)
                \cos\left(\theta\,\rho\frac{\partial}{\partial\rho}\right)f(\rho) &=& \cos(Y\partial_X)f(X), \\
                \sin\left(\theta\,\rho\frac{\partial}{\partial\rho}\right)f(\rho) &=&
                \sin(Y\partial_X)f(X).
\end{eqnarray*}}\hfill\parbox{1cm}{\begin{eqnarray}\end{eqnarray}}

\textbf{Proof.} Let
$a_1=e^{i\theta_1},\,\ldots,\,a_n=e^{i\theta_n}$, according to
Theorem 2.4, then
\[\exp\left(i\theta\,\rho\frac{\partial}{\partial\rho}\right)f(\rho)=f(\rho{e}^{i\theta})=f(X+iY)=\exp(iY\partial_X)f(X),\qquad X\in\mathbb{R}^n,\;\;Y\in\mathbb{R}_n.\]
Considering
\begin{eqnarray*}
% \nonumber to remove numbering (before each equation)
  \exp\left(i\theta\,\rho\frac{\partial}{\partial\rho}\right) &=& \cos\left(\theta\,\rho\frac{\partial}{\partial\rho}\right)
  +i\sin\left(\theta\,\rho\frac{\partial}{\partial\rho}\right), \\
  \exp(iY\partial_X) &=& \cos(Y\partial_X)+i\sin(Y\partial_X),
\end{eqnarray*}
thus we have Theorem 2.5, which transforms the abstract operators
taking $\rho\partial_\rho$ as the operator element
into those taking $\partial_X$ as the operator element.

It is easy to see that we have following Corollary:

\textbf{Corollary 2.6.} Let $f(z)\in{C^\infty}(\Omega),\,z=x+iy\in\Omega\subset\mathbb{C}^1$ be an
arbitrary analytic function, then the harmonic functions $u(x,y)$ and $\upsilon(x,y)$ on complex plane can be expressed as

\parbox{12.5cm}{\begin{eqnarray*}\label{15**}
u(x,y)=\cos\left(y\frac{\partial}{\partial{x}}\right)f(x)\!\!
&=& \!\!\cos\left(\theta\rho\frac{\partial}{\partial\rho}\right)f(\rho),\\
\upsilon(x,y)=\sin\left(y\frac{\partial}{\partial{x}}\right)f(x)\!\!
&=& \!\!\sin\left(\theta\rho\frac{\partial}{\partial\rho}\right)f(\rho),
\end{eqnarray*}}\hfill\parbox{1cm}{\begin{eqnarray}\end{eqnarray}}
where $\rho=\sqrt{x^2+y^2},\,\theta=\arctan\frac{y}{x}$. Especially, if $y=kx$, then (\ref{15**}) gives

\parbox{13cm}{\begin{eqnarray*}\label{15*}
\left.\cos\left(y\frac{\partial}{\partial{x}}\right)f(x)\right|_{y=kx}
&=&
 \beta^{x\frac{\partial}{\partial{x}}}\cos\left(\alpha{x}\frac{\partial}{\partial{x}}\right)f(x),\\
\left.\sin\left(y\frac{\partial}{\partial{x}}\right)f(x)\right|_{y=kx}
&=&
\beta^{x\frac{\partial}{\partial{x}}}\sin\left(\alpha{x}\frac{\partial}{\partial{x}}\right)f(x),
\end{eqnarray*}}\hfill\parbox{1cm}{\begin{eqnarray}\end{eqnarray}}
where $\beta=\sqrt{1+k^2},\;\alpha=\arctan{k}$, which can be used to solve the boundary value problems of 2-dimensional Laplace equation on polygonal domains.

\subsection{Integral representation of abstract operators}
\noindent

Clearly, such algorithms as (\ref{y1})-(\ref{ey3}) do not exist for
most abstract operators. Therefore, for more complex abstract operators, establishing the integral expression is important way to extend algorithms of the abstract operators. Guang-Qing Bi has obtained the following results in 1997:

\textbf{Theorem BI3.} (See \cite{bi97},\cite{bi99}) Let $P(\partial_x)$ be an $m$-order partial
differential operators of any kind and there exist
$a_1,\,a_2,\,\cdots,\,a_k$ of real and partial differential
operators $A_1,\,A_2,\,\cdots,\,A_k$ of the order less than
$[(m+1)/2]$ such that
$P(\partial_x)\equiv{a}_1A_1^2+a_2A_2^2+\cdots+a_kA_k^2$. If
$k=2\nu+3,\;\nu=0,\,1,\,2,\cdots$, then
$\forall{f(x)}\in{C^\infty}(\Omega),\,x\in\Omega\subset\mathbb{R}^n$ we have

\begin{eqnarray}\label{y6}
% \nonumber to remove numbering (before each equation)
\frac{\sinh\left(tP(\partial_x)^{1/2}\right)}{P(\partial_x)^{1/2}}f(x)
  &=& t\underbrace{\int^t_0tdt\cdots}_\nu\int^t_0\,tdt\frac{(P(\partial_x))^{\nu}}{2^{\nu+2}\pi^{\nu+1}}\nonumber\\
   & & \times\,\int^\pi_{-\pi}\underbrace{\int^\pi_0\cdots}_{k-2}\int^\pi_0e^{\eta_1a_1^{1/2}A_1+\cdots+\eta_ka_k^{1/2}A_k}f(x)\,d\sigma_k\nonumber\\
   & & +\,\sum^{\nu-1}_{i=0}\frac{t^{2i+1}}{(2i+1)!}(P(\partial_x))^{i}f(x).
\end{eqnarray}
Similarly
\begin{eqnarray}\label{y6'}
% \nonumber to remove numbering (before each equation)
\frac{\sin\left({t}P(\partial_x)^{1/2}\right)}{P(\partial_x)^{1/2}}f(x)
   &=& t\underbrace{\int^t_0tdt\cdots}_\nu\int^t_0tdt[-P(\partial_x)]^{\nu}\nonumber\\
& & \times\,\left(\frac{2}{\pi}\right)^{\nu+1}\underbrace{\int^{\pi/2}_0\cdots}_{k-1}\int^{\pi/2}_0\cos(\eta_1a_1^{1/2}A_1)\cdots\cos(\eta_ka_k^{1/2}A_k)f(x)\,d\sigma_k\nonumber\\
   & & +\,\sum^{\nu-1}_{i=0}\frac{t^{2i+1}}{(2i+1)!}[-P(\partial_x)]^{i}f(x).
\end{eqnarray}
Where $x\in\mathbb{R}^n,\;t\in\mathbb{R}^1.\;\;\eta\in\mathbb{R}_k$ is the integral variable and
\begin{eqnarray*}
% \nonumber to remove numbering (before each equation)
  \eta_1 &=& t\cos\theta_1,\\
  \eta_2 &=& t\sin\theta_1\cos\theta_2,\\
  \eta_3 &=& t\sin\theta_1\sin\theta_2\cos\theta_3,\\
         &\cdots&\\
  \eta_p &=& t\sin\theta_1\sin\theta_2\cdots\sin\theta_{p-1}\cos\theta_p,\\
  \eta_{p+1} &=& t\sin\theta_1\sin\theta_2\cdots\sin\theta_p\cos\phi,\\
  \eta_{p+2} &=&
  \eta_k\;=\;t\sin\theta_1\sin\theta_2\cdots\sin\theta_p\sin\phi;
\end{eqnarray*}
\[d\sigma_k=\sin^{k-2}\theta_1\sin^{k-3}\theta_2\cdots\sin\theta_{k-2}d\theta_1d\theta_2\cdots{d}\theta_{k-2}d\phi.\]

\textbf{Proof.} In (\ref{y6}), let
$f(x)=e^{\xi{x}},\,x\in\mathbb{R}^n,\xi\in\mathbb{R}_n$, and the
symbols of the partial differential operators $A_j,j=1,2,\cdots,k$ be
denoted by $\chi_j(\xi),\;\beta_j=a_j^{1/2}\chi_j(\xi)$, then (\ref{y6}) degenerates to its characteristic
equation by Definition 2.4, namely
\begin{eqnarray}\label{y6z}
% \nonumber to remove numbering (before each equation)
\frac{\sinh\left(tP(\xi)^{1/2}\right)}{P(\xi)^{1/2}}&=&
t\underbrace{\int^t_0tdt\cdots}_\nu\int^t_0\,tdt\frac{(P(\xi))^{\nu}}{2^{\nu+2}\pi^{\nu+1}}
\int^\pi_{-\pi}\underbrace{\int^\pi_0\cdots}_{k-2}\int^\pi_0e^{\eta_1\beta_1+\cdots+\eta_k\beta_k}d\sigma_k\nonumber\\
   & &
   +\,\sum^{\nu-1}_{i=0}\frac{t^{2i+1}}{(2i+1)!}(P(\xi))^i,\quad\nu=\frac{k-3}{2}.
\end{eqnarray}
According to the Analytic continuous fundamental theorem, we only
need to prove (\ref{y6z}). Solving the integral on a hypersphere on
the right side of (\ref{y6z}), we have
\begin{eqnarray*}
\frac{\sinh\left(tP(\xi)^{1/2}\right)}{P(\xi)^{1/2}}=
t\underbrace{\int^t_0tdt\cdots}_\nu\int^t_0tdt
\sum^\infty_{j=0}\frac{(P(\xi))^{\nu+j}t^{2j}}{(2j)!!(2j+2\nu+1)!!}+\sum^{\nu-1}_{i=0}\frac{t^{2i+1}(P(\xi))^i}{(2i+1)!}.
\end{eqnarray*}
Then it is proved by the termwise integration of the infinite series
on the right side of the equality. Similarly, we have (\ref{y6'}).

When $A_1,\,\cdots,\,A_k$ in the right-hand of (\ref{y6}) and (\ref{y6'}) are one
order partial differential operators, then the abstract operators
$e^{\eta_1a_1^{1/2}A_1+\cdots+\eta_ka_k^{1/2}A_k}$ and $\cos(\eta_1a_1^{1/2}A_1)\cdots\cos(\eta_ka_k^{1/2}A_k)$ are one of the
following five simplest operators:
\[\exp(h\partial_x),\;\sin(h\partial_x),\;\cos(h\partial_x),\;\sinh(h\partial_x)\;\mbox{and}\;\cosh(h\partial_x).\]

According to (\ref{y6}), we can easily derive:

\textbf{Example 2.7.} Let $\Delta_n$ be the n-dimensional Laplacian, $n-2=2\nu+1,\,a\in\mathbb{R},\,t\in\mathbb{R}_+^1$,
then we have
\begin{eqnarray}\label{37}
% \nonumber to remove numbering (before each equation)
\frac{\sinh\left(at\Delta_n{}^{1/2}\right)}{a\Delta_n{}^{1/2}}f(x)&=&
t\underbrace{\int^t_0tdt\cdots}_\nu\int^t_0\frac{(a^{2}\Delta_n)^{\nu}}{S_n}\int_{S_n}f(\xi)\,dS_n\,tdt\nonumber\\
   & &
   +\,\sum^{\nu-1}_{i=0}\frac{t^{2i+1}}{(2i+1)!}(a^{2}\Delta_n)^{i}f(x),\quad\forall{f(x)}\in{C}^{2\nu}(\Omega).
\end{eqnarray}
Where $S_n=2(2\pi)^{\nu+1}(at)^{n-1}$, $\xi\in\mathbb{R}_n$ is the
integral variable on the hypersphere
$(\xi_1-x_1)^2+(\xi_2-x_2)^2+\cdots+(\xi_n-x_n)^2=(at)^2$, and
$dS_n$ is its surface element.

If $A_1,\,\cdots,\,A_k$ in (\ref{y6}) and (\ref{y6'}) are partial
differential operators of the order great than 1, then the order can
be lowered by taking the following Theorem:

\textbf{Theorem BI4.} (See \cite{bi97}) Let $P(\partial_x),\;x\in\mathbb{R}^n$ be a
partial differential operator of any order and
$f(x)\in{C}^\infty(\Omega),\,x\in\Omega\subset\mathbb{R}^n$ which make the integral in the follow
formula meaningful, then
\begin{equation}\label{y7}
    e^{\lambda
    P(\partial_x)}f(x)=\frac{1}{2\sqrt{\pi}}\int^\infty_{-\infty}e^{-\zeta^2/4}e^{\lambda^{1/2}\zeta
    P(\partial_x)^{1/2}}f(x)d\zeta,\quad\forall\lambda\in\mathbb{C}^1.
\end{equation}

According to (\ref{y7}) and Corollary 2.1, we can easily derive:

\textbf{Example 2.8.} Let
$h_{\lambda,a}(\zeta)=a\lambda+\sqrt{\lambda/2}\zeta,\;\forall\lambda,a\in\mathbb{C},\;x\in\mathbb{R}^n,\;t\in\mathbb{R}^1$,
then we have
\[\exp\left(-a^2t\frac{\partial^2}{\partial{x}_j^2}\right)g(x)=
\frac{1}{2\sqrt{\pi}}\int^\infty_{-\infty}e^{-\zeta^2/4}\cos\left(a\sqrt{t}\,\zeta\frac{\partial}{\partial{x_j}}\right)g(x)d\zeta.\]
\[\cos\left(\lambda\frac{\partial^2}{\partial{x}_j^2}\right)g(x)=\frac{1}{2\sqrt{\pi}}\int^\infty_{-\infty}e^{-\zeta^2/4}
e^{\sqrt{\lambda/2}\,\zeta\frac{\partial}{\partial{x}_j}}\cos\left(\sqrt{\frac{\lambda}{2}}\,\zeta\frac{\partial}{\partial{x}_j}\right)g(x)d\zeta.\]
\[\sin\left(a\lambda\frac{\partial}{\partial{x}_j}+\lambda\frac{\partial^2}{\partial{x}_j^2}\right)g(x)=\frac{1}{2\sqrt{\pi}}\int^\infty_{-\infty}e^{-\zeta^2/4}
e^{\sqrt{\lambda/2}\,\zeta\frac{\partial}{\partial{x}_j}}\sin\left(h_{\lambda,a}(\zeta)\frac{\partial}{\partial{x}_j}\right)g(x)d\zeta.\]

\section{Applications of abstract operators to partial differential equations}
\noindent

Solving the ordinary or partial differential equations, is
constructing the algorithms of the inverse operators of ordinary or
partial differential operators.

\subsection{Laplace transform of $n+1$ dimensional partial differential equations}
\noindent

In terms of abstract operators, solving the initial value problem of $n+1$ dimensional
partial differential equations is similar to solving the ordinary differential equations with respect to the
variable $t$, thus we can introduce the Laplace transform to further
simplify the solving process. Now let us prove that if we combine abstract operators
with the Laplace transform, then no matter how complex the initial value problem of linear partial differential
equations is, the solving process is simply easy.

\textbf{Theorem 3.1.} Puppose that $f^{(k)}\!\left(\lambda\frac{\partial}{\partial{s}}\right),\,k=0,1,2,\cdots$ are the abstract operators and $s,\lambda$ are the complex parameter, if
$\exists{v(s),u(s)}\in{C^\infty}(\Omega),\;s\in\Omega\subset\mathbb{C}^1$ can
make the infinite series on the right side of (\ref{t2}) uniform
convergent, then it will uniform converges to the left side of this
equality, namely
\begin{equation}\label{t2}
f\!\left(\lambda\frac{\partial}{\partial{s}}\right)(vu)=
\sum_{k=0}^\infty\frac{\lambda^k}{k!}\,\frac{\partial^kv}{\partial{s^k}}\,f^{(k)}\!\left(\lambda\frac{\partial}{\partial{s}}\right)u,\quad|\lambda|<R.
\end{equation}

\textbf{Proof.} Taking Taylor formula (\ref{2})
as its characteristic equation, and using the Analytic continuous fundamental theorem, then (\ref{t2}) can be easily proved.

\textbf{Corollary 3.1.} In (\ref{t2}), when
$\lambda=-1,\;v=s,\;u=1/s,\,s\neq0$, we have
\begin{equation}\label{t3}
f'\left(-\frac{\partial}{\partial{s}}\right)\frac{1}{s}=sf\left(-\frac{\partial}{\partial{s}}\right)\frac{1}{s}-f(0)
\quad\mbox{or}\quad\mathscr{L}f'(t)=s\mathscr{L}f(t)-f(0).
\end{equation}
Where $\mathscr{L}f(t)$ denotes
\begin{equation}\label{t4}
\mathscr{L}f(t)=f\left(-\frac{\partial}{\partial{s}}\right)\frac{1}{s},\quad\;\Re(s)>0,\;t\in\mathbb{R}_+^1=\{t\in\mathbb{R}^1|t>0\}.
\end{equation}
This is the internal connection between abstract operators and
Laplace transform, thus the abstract operators becomes a significant
tool to compute Laplace transform. For instance, by (\ref{t4}) we
have
\begin{equation}\label{t5}
\mathscr{L}[g(t)f(t)]=g\left(-\frac{\partial}{\partial{s}}\right)F(s),\quad(F(s)=\mathscr{L}f(t)).
\end{equation}

So the symbols of the abstract operators on the right side of (\ref{t4}) and (\ref{t5}) can be further
described by using conditions of the Laplace transform.

\textbf{Theorem 3.2.} Let $m\geq1$,
$P(\partial_x)$ be an arbitrary order partial differential equations, then we have
\begin{equation}\label{32}
    \left\{\begin{array}{l@{\qquad}l}\displaystyle
    \left(\frac{\partial^2}{\partial{t^2}}-P(\partial_x)\right)^mu=f(x,t),&x\in\mathbb{R}^n,\;t\in\mathbb{R}_+^1,\\\displaystyle
    \left.\frac{\partial^ju}{\partial{t^j}}\right|_{t=0}=\varphi_j(x),&j=0,1,2,\ldots,2m-1.
    \end{array}\right.
\end{equation}
\begin{eqnarray}\label{32'}
  u(x,t) &=&
  \int^t_0\int^{t-\tau}_0\frac{\left[(t-\tau)^2-\tau'^2\right]^{m-2}}{(2m-2)!!\,(2m-4)!!}\,
  \frac{\sinh\left(\tau'P(\partial_x)^{1/2}\right)}{P(\partial_x)^{1/2}}\,f(x,\tau)\,\tau'd\tau'\,d\tau \nonumber\\
   & & +\,\sum^{m-1}_{k=0}(-1)^k{m\choose{k}}P(\partial_x)^k\sum^{2m-1-2k}_{j=0}\frac{\partial^{2m-1-2k-j}}{\partial t^{2m-1-2k-j}}\int^t_0
   \frac{(t^2-\tau^2)^{m-2}\tau}{(2m-2)!!\,(2m-4)!!}\nonumber\\
   & & \times\,\frac{\sinh\left(\tau P(\partial_x)^{1/2}\right)}{P(\partial_x)^{1/2}}\,\varphi_j(x)\,d\tau.
\end{eqnarray}

\textbf{Proof.} The Laplace transform of the partial differential
equations, with respect to $t$ and considering the initial condition,
is
\[\sum^m_{k=0}(-1)^k{m\choose{k}}P(\partial_x)^k\left(s^{2m-2k}U(x,s)-\sum^{2m-1-2k}_{j=0}s^{2m-1-2k-j}\varphi_j(x)\right)=F(x,s).\]
Where $U(x,s)=\mathscr{L}u(x,t),\;F(x,s)=\mathscr{L}f(x,t)$. Let
$G_m(\partial_x,t)=\mathscr{L}^{-1}[1/(s^2-P(\partial_x))^m]$,
solving $U(x,s)$ and its inverse Laplace transform is
\begin{eqnarray*}
% \nonumber to remove numbering (before each equation)
   u(x,t)&=& \mathscr{L}^{-1}U(x,s)= \mathscr{L}^{-1}\frac{1}{(s^2-P(\partial_x))^m}F(x,s) \\
   & &+\,\mathscr{L}^{-1}\sum^{m-1}_{k=0}(-1)^k{m\choose{k}}P(\partial_x)^k\sum^{2m-1-2k}_{j=0}\frac{s^{2m-1-2k-j}}{(s^2-P(\partial_x))^m}\varphi_j(x) \\
   &=& G_m(\partial_x,t)\ast f(x,t)\\
   & &+\,\sum^{m-1}_{k=0}(-1)^k{m\choose{k}}P(\partial_x)^k
   \sum^{2m-1-2k}_{j=0}\frac{\partial^{2m-1-2k-j}}{\partial
   t^{2m-1-2k-j}}G_m(\partial_x,t)\varphi_j(x).
\end{eqnarray*}

Now let us solve $G_m(\partial_x,t)$.
Puppose that $f(t)$ is a real- or complex-valued function of the (time) variable $t>0$ and $s$ is a real or complex parameter, we have
\[\mathscr{L}\,\left(\int^t_0\!\cdot\,{t}dt\right)^{m-1}\!\!f(t)=\left(-\frac{1}{s}\frac{\partial}{\partial{s}}\right)^{m-1}\mathscr{L}f(t),\quad\forall{m}\geq1.\]

Let $f(t)=\sin{bt},\,b\in\mathbb{C},\;t\in\mathbb{R}_+$, we
have
\[\left(\int^t_0\cdot\,tdt\right)^{m-1}\!\!\sin{bt}=\mathscr{L}^{-1}\left(-\frac{1}{s}\frac{\partial}{\partial{s}}\right)^{m-1}\!\!\frac{b}{s^2+b^2}
=\mathscr{L}^{-1}\frac{2^{m-1}(m-1)!}{(s^2+b^2)^m}b.\] Let
$b=iP(\xi)^{1/2},\;\xi\in\mathbb{R}_n$, we have
\[\mathscr{L}^{-1}\frac{1}{(s^2-P(\xi))^m}=\frac{1}{(2m-2)!!}\left(\int^t_0\cdot\,tdt\right)^{m-1}\frac{\sinh\left(tP(\xi)^{1/2}\right)}{P(\xi)^{1/2}}.\]
Taking this one as the characteristic  equation, according to the
Analytic continuous fundamental theorem, we have
\begin{equation}\label{33}
G_m(\partial_x,t)=\mathscr{L}^{-1}\frac{1}{(s^2-P(\partial_x))^m}=
\frac{1}{(2m-2)!!}\left(\int^t_0\cdot\,tdt\right)^{m-1}\frac{\sinh\left(tP(\partial_x)^{1/2}\right)}{P(\partial_x)^{1/2}}.
\end{equation}
By Reference \cite{bi01}, we can easily derive the
following integral formula
\begin{equation}\label{34}
\left(\int^t_a\cdot\,tdt\right)^{m}f(t)=\underbrace{\int^t_atdt\cdots}_{m}\int^t_af(t)\,tdt=\int^t_a\frac{(t^2-\tau^2)^{m-1}}{(2m-2)!!}f(\tau)\,\tau d\tau.
\end{equation}
Applying (\ref{34}) to (\ref{33}), we have the expression of abstract operators $G_m(\partial_x,t)$:
\begin{equation}\label{35}
G_m(\partial_x,t)g(x)=
\int^t_0\frac{(t^2-\tau'^2)^{m-2}}{(2m-2)!!\,(2m-4)!!}\,\frac{\sinh\left(\tau'P(\partial_x)^{1/2}\right)}{P(\partial_x)^{1/2}}g(x)\,\tau'd\tau'.
\end{equation}
Thus Theorem 3.2 is proved.

It is extremely complex to prove Theorem 3.2 even if $\varphi_j(x)=0,\,j=0,1,\ldots,2m-1$ without using the Laplace transform (See \cite{bi01}).

In 1999, Guang-Qing Bi has obtained the following results (See Reference \cite{bi99}):

\textbf{Theorem BI5.} (See \cite{bi99}) Let $a_1,a_2,\ldots,a_m$ be arbitrary real or
complex numbers different from each other, $P(\partial_x)$ be a
partial differential operators of any order, then we have
\begin{equation}\label{I1}
    \left\{\begin{array}{l@{\qquad}l}\displaystyle
    \prod^m_{i=1}(\frac{\partial}{\partial{t}}-a_iP(\partial_x))u=f(x,t),&x\in\mathbb{R}^n,\;t\in\mathbb{R}_+^1,\;\forall{m}\geq1,\\\displaystyle
    \left.\frac{\partial^ju}{\partial{t^j}}\right|_{t=0}=0,&j=0,1,2,\ldots,m-1.
    \end{array}\right.
\end{equation}
\begin{equation}\label{I1'}
u(x,t)=\int^t_0\int^{t-\tau}_0\frac{(t-\tau-\tau')^{m-2}}{(m-2)!}\sum^m_{j=1}\frac{a_j^{m-1}}{\prod^m_{i=1\atop
i\neq{j}}(a_j-a_i)} e^{\tau'a_jP(\partial_x)}f(x,\tau)\,d\tau'd\tau.
\end{equation}

\textbf{Theorem BI6.} (See \cite{bi99}) Let $a_1,a_2,\ldots,a_m$ be arbitrary real or
complex numbers different from each other, $P(\partial_x)$ be a
partial differential operators of any order, then we have
\begin{equation}\label{I2}
    \left\{\begin{array}{l@{\qquad}l}\displaystyle
    \prod^m_{i=1}(\frac{\partial^2}{\partial{t^2}}-a_i^2P(\partial_x))u=f(x,t),&x\in\mathbb{R}^n,\;t\in\mathbb{R}_+^1,\;\forall{m}\geq1,\\\displaystyle
    \left.\frac{\partial^ju}{\partial{t^j}}\right|_{t=0}=0,&j=0,1,2,\ldots,2m-1.
    \end{array}\right.
\end{equation}
\begin{equation}\label{I2'}
u(x,t)=\int^t_0\int^{t-\tau}_0\frac{(t-\tau-\tau')^{2m-3}}{(2m-3)!}\sum^m_{j=1}\frac{a_j^{2m-2}}{\prod^m_{i=1
\atop i\neq{j}}(a_j^2-a_i^2)}
\frac{\sinh(\tau'a_jP(\partial_x)^{1/2})}{a_jP(\partial_x)^{1/2}}f(x,\tau)\,d\tau'd\tau.
\end{equation}

On this basis, by using the abstract operators and Laplace transform we have obtained the following Theorems:

\textbf{Theorem 3.3.} Let $a_1,a_2,\ldots,a_m$ be arbitrary real or
complex roots different from each other for
$b_0+b_1\chi+b_2\chi^2+\cdots+b_m\chi^m=0$, and
$P(\partial_x,\partial_t)$ be the partial differential operators
defined by
$$P(\partial_x,\partial_t)=\sum^m_{k=0}b_kP(\partial_x)^{m-k}\frac{\partial^k}{\partial{t^k}},
\quad{x}\in\mathbb{R}^n,\;t\in\mathbb{R}_+^1,\;\forall{m}\geq1.$$
Where $P(\partial_x)$ is a partial differential operators of any
order. Then we have
\begin{equation}\label{j1}
    \left\{\begin{array}{l@{\qquad}l}\displaystyle
    P(\partial_x,\partial_t)u=f(x,t),&x\in\mathbb{R}^n,\;t\in\mathbb{R}_+^1,\\\displaystyle
    \left.\frac{\partial^ru}{\partial{t^r}}\right|_{t=0}=\varphi_r(x),&r=0,1,2,\ldots,m-1.
    \end{array}\right.
\end{equation}
\begin{eqnarray}\label{j1'}
  u(x,t) &=&
  \int^t_0\int^{t-\tau}_0\frac{(t-\tau-\tau')^{m-2}}{(m-2)!}\sum^m_{j=1}\frac{a_j^{m-1}}{\prod^m_{i=1 \atop i\neq{j}}(a_j-a_i)}\,
e^{\tau'a_jP(\partial_x)}f(x,\tau)\,d\tau'd\tau\nonumber\\
   & & +\,\sum^m_{k=1}b_kP(\partial_x)^{m-k}\sum^{k-1}_{r=0}\frac{\partial^{k-1-r}}{\partial t^{k-1-r}}\int^t_0
   \frac{(t-\tau)^{m-2}}{(m-2)!}\nonumber\\
   & & \times\,\sum^m_{j=1}\frac{a_j^{m-1}}{\prod^m_{i=1 \atop i\neq{j}}(a_j-a_i)}\,e^{\tau{a_j}P(\partial_x)}\,\varphi_r(x)\,d\tau.
\end{eqnarray}

\textbf{Theorem 3.4.} Let $a_1,a_2,\ldots,a_m$ be arbitrary real or
complex roots different from each other, satisfy
$\sum^m_{k=0}b_{2k}\chi^{2k}=\prod^m_{i=1}(\chi^2-a_i^2)$, and
$P(\partial_x,\partial_t)$ be the partial differential operators
defined by
$$P(\partial_x,\partial_t)=\sum^m_{k=0}b_{2k}P(\partial_x)^{m-k}\frac{\partial^{2k}}{\partial{t^{2k}}},
\quad{x}\in\mathbb{R}^n,\;t\in\mathbb{R}_+^1,\;\forall{m}\geq1.$$
Where $P(\partial_x)$ be a partial differential operators of any
order. Then we have
\begin{equation}\label{j2}
    \left\{\begin{array}{l@{\qquad}l}\displaystyle
    P(\partial_x,\partial_t)u=f(x,t),&x\in\mathbb{R}^n,\;t\in\mathbb{R}_+^1,\\\displaystyle
    \left.\frac{\partial^ru}{\partial{t^r}}\right|_{t=0}=\varphi_r(x),&r=0,1,2,\ldots,2m-1.
    \end{array}\right.
\end{equation}
\begin{eqnarray}\label{j2'}
  u(x,t) &=&
  \int^t_0\int^{t-\tau}_0\frac{(t-\tau-\tau')^{2m-3}}{(2m-3)!}\sum^m_{j=1}\frac{a_j^{2m-2}}{\prod^m_{i=1 \atop i\neq{j}}(a_j^2-a_i^2)}
\frac{\sinh(\tau'a_jP(\partial_x)^{1/2})}{a_jP(\partial_x)^{1/2}}f(x,\tau)\,d\tau'd\tau\nonumber\\
   & & +\,\sum^m_{k=1}b_{2k}P(\partial_x)^{m-k}\sum^{2k-1}_{r=0}\frac{\partial^{2k-1-r}}{\partial t^{2k-1-r}}\int^t_0
   \frac{(t-\tau)^{2m-3}}{(2m-3)!}\nonumber\\
   & & \times\,\sum^m_{j=1}\frac{a_j^{2m-2}}{\prod^m_{i=1 \atop i\neq{j}}(a_j^2-a_i^2)}
   \frac{\sinh(\tau{a_j}P(\partial_x)^{1/2})}{a_jP(\partial_x)^{1/2}}\,\varphi_r(x)\,d\tau.
\end{eqnarray}

Let us prove these two Theorems. According to the Theorem BI5 and
Theorem BI6, we just need to prove the following corollary of the
Theorem 3.3 and Theorem 3.4:

\textbf{Corollary 3.2.} Let $a_1,a_2,\ldots,a_m$ be arbitrary real or
complex roots different from each other for
$b_0+b_1\chi+b_2\chi^2+\cdots+b_m\chi^m=0$, and
$P(\partial_x,\partial_t)$ be the partial differential operators
defined by
$$P(\partial_x,\partial_t)=\sum^m_{k=0}b_kP(\partial_x)^{m-k}\frac{\partial^k}{\partial{t^k}},
\quad{x}\in\mathbb{R}^n,\;t\in\mathbb{R}_+^1,\;\forall{m}\geq1.$$
Where $P(\partial_x)$ be a partial differential operators of any
order. Then we have
\begin{equation}\label{j3}
    \left\{\begin{array}{l@{\qquad}l}\displaystyle
    P(\partial_x,\partial_t)u=0,&x\in\mathbb{R}^n,\;t\in\mathbb{R}_+^1,\\\displaystyle
    \left.\frac{\partial^ru}{\partial{t^r}}\right|_{t=0}=\varphi_r(x),&r=0,1,2,\ldots,m-1.
    \end{array}\right.
\end{equation}
\begin{eqnarray}\label{j3'}
  u(x,t)=\sum^m_{k=1}b_kP(\partial_x)^{m-k}\sum^{k-1}_{r=0}\frac{\partial^{k-1-r}}{\partial
t^{k-1-r}}\int^t_0\frac{(t-\tau)^{m-2}}{(m-2)!}
\sum^m_{j=1}\frac{a_j^{m-1}e^{\tau{a_j}P(\partial_x)}}{\prod^m_{i=1
\atop i\neq{j}}(a_j-a_i)}\,\varphi_r(x)\,d\tau.
\end{eqnarray}

\textbf{Proof.} Considering initial conditions, the Laplace
transform of the Eq (\ref{j3}) with respect to $t$ is
$$\sum^m_{k=0}b_kP(\partial_x)^{m-k}\left(s^kU(x,s)-\sum^{k-1}_{r=0}s^{k-1-r}\varphi_r(x)\right)=0.$$
Where $U(x,s)=\mathscr{L}u(x,t)$, considering
$\prod^m_{i=1}(s-a_iP(\partial_x))=\sum^m_{k=0}b_ks^kP(\partial_x)^{m-k}$
we have
$$\prod^m_{i=1}(s-a_iP(\partial_x))U(x,s)-\sum^m_{k=1}b_kP(\partial_x)^{m-k}\sum^{k-1}_{r=0}s^{k-1-r}\varphi_r(x)=0.$$
We need to introduce an abstract operators $G_m(\partial_x,t)$,
defined as
$$G_m(\partial_x,t)=\mathscr{L}^{-1}\frac{1}{\prod^m_{i=1}(s-a_iP(\partial_x))}.$$
By solving $U(x,s)$, we have its inverse Laplace transform:
\begin{eqnarray}\label{j3''}
% \nonumber to remove numbering (before each equation)
u(x,t)&=& \mathscr{L}^{-1}U(x,s)= \sum^m_{k=1}b_kP(\partial_x)^{m-k}\sum^{k-1}_{r=0}\mathscr{L}^{-1}\frac{s^{k-1-r}}{\prod^m_{i=1}(s-a_iP(\partial_x))}\,\varphi_r(x)\nonumber\\
   &=& \sum^m_{k=1}b_kP(\partial_x)^{m-k}\sum^{k-1}_{r=0}\frac{\partial^{k-1-r}}{\partial{t^{k-1-r}}}\,G_m(\partial_x,t)\varphi_r(x).
\end{eqnarray}

Now let us solve $G_m(\partial_x,t)$. Considering initial
conditions, the Laplace transform of the Eq (\ref{I1}) with respect
to $t$ is
$$\prod^m_{i=1}(s-a_iP(\partial_x))U(x,s)=F(x,s),$$
where $F(x,s)=\mathscr{L}f(x,t)$. By solving $U(x,s)$ and using the
convolution theorem, we have its inverse Laplace transform:
$$u(x,t)=\mathscr{L}^{-1}U(x,s)=\mathscr{L}^{-1}\frac{1}{\prod^m_{i=1}(s-a_iP(\partial_x))}F(x,s)=G_m(\partial_x,t)*f(x,t).$$
By comparing (\ref{I1'}) with $u(x,t)=G_m(\partial_x,t)*f(x,t)$, we
have the expression of the abstract operators $G_m(\partial_x,t)$:
\begin{equation}\label{j4}
G_m(\partial_x,t)=\int^t_0\frac{(t-\tau)^{m-2}}{(m-2)!}
\sum^m_{j=1}\frac{a_j^{m-1}}{\prod^m_{i=1 \atop
i\neq{j}}(a_j-a_i)}e^{\tau{a_j}P(\partial_x)}d\tau.
\end{equation}
Applying (\ref{j4}) to (\ref{j3''}), thus Corollary 3.2 is proved.

\textbf{Corollary 3.3.} Let $a_1,a_2,\ldots,a_m$ be arbitrary real or
complex roots different from each other, which satisfy
$\sum^m_{k=0}b_{2k}\chi^{2k}=\prod^m_{i=1}(\chi^2-a_i^2)$, and
$P(\partial_x,\partial_t)$ be the partial differential operators
defined by
$$P(\partial_x,\partial_t)=\sum^m_{k=0}b_{2k}P(\partial_x)^{m-k}\frac{\partial^{2k}}{\partial{t^{2k}}},
\quad{x}\in\mathbb{R}^n,\;t\in\mathbb{R}_+^1,\;\forall{m}\geq1.$$
Where $P(\partial_x)$ is a partial differential operators of any
order. Then we have
\begin{equation}\label{j5}
    \left\{\begin{array}{l@{\qquad}l}\displaystyle
    P(\partial_x,\partial_t)u=0,&x\in\mathbb{R}^n,\;t\in\mathbb{R}_+^1,\\\displaystyle
    \left.\frac{\partial^ru}{\partial{t^r}}\right|_{t=0}=\varphi_r(x),&r=0,1,2,\ldots,2m-1.
    \end{array}\right.
\end{equation}
\begin{eqnarray}\label{j5'}
  u(x,t) &=&
   \sum^m_{k=1}b_{2k}P(\partial_x)^{m-k}\sum^{2k-1}_{r=0}\frac{\partial^{2k-1-r}}{\partial t^{2k-1-r}}\int^t_0
   \frac{(t-\tau)^{2m-3}}{(2m-3)!}\nonumber\\
   & & \times\,\sum^m_{j=1}\frac{a_j^{2m-2}}{\prod^m_{i=1 \atop i\neq{j}}(a_j^2-a_i^2)}
   \frac{\sinh(\tau{a_j}P(\partial_x)^{1/2})}{a_jP(\partial_x)^{1/2}}\,\varphi_r(x)\,d\tau.
\end{eqnarray}

\textbf{Proof.} Considering initial conditions, the Laplace
transform of the Eq (\ref{j5}) with respect to $t$ is
$$\sum^m_{k=0}b_{2k}P(\partial_x)^{m-k}\left(s^{2k}U(x,s)-\sum^{2k-1}_{r=0}s^{2k-1-r}\varphi_r(x)\right)=0.$$
Where $U(x,s)=\mathscr{L}u(x,t)$, considering
$$\prod^m_{i=1}(s^2-a_i^2P(\partial_x))=\sum^m_{k=0}b_{2k}s^{2k}P(\partial_x)^{m-k}$$
we have
$$\prod^m_{i=1}(s^2-a_i^2P(\partial_x))U(x,s)-\sum^m_{k=1}b_{2k}P(\partial_x)^{m-k}\sum^{2k-1}_{r=0}s^{2k-1-r}\varphi_r(x)=0.$$
We need to introduce an abstract operators $G_m(\partial_x,t)$,
defined as
$$G_m(\partial_x,t)=\mathscr{L}^{-1}\frac{1}{\prod^m_{i=1}(s^2-a_i^2P(\partial_x))}.$$
By solving $U(x,s)$, we have its inverse Laplace transform:
\begin{eqnarray}\label{j5''}
% \nonumber to remove numbering (before each equation)
   u(x,t)&=& \mathscr{L}^{-1}U(x,s)= \sum^m_{k=1}b_{2k}P(\partial_x)^{m-k}\sum^{2k-1}_{r=0}\mathscr{L}^{-1}\frac{s^{2k-1-r}}{\prod^m_{i=1}(s^2-a_i^2P(\partial_x))}\,
   \varphi_r(x)\nonumber\\
   &=& \sum^m_{k=1}b_{2k}P(\partial_x)^{m-k}\sum^{2k-1}_{r=0}\frac{\partial^{2k-1-r}}{\partial{t^{2k-1-r}}}\,G_m(\partial_x,t)\varphi_r(x).
\end{eqnarray}

Now let us solve the $G_m(\partial_x,t)$. Considering initial
conditions, the Laplace transform of the Eq (\ref{I2}) with respect
to $t$ is
$$\prod^m_{i=1}(s^2-a_i^2P(\partial_x))U(x,s)=F(x,s),$$
where $F(x,s)=\mathscr{L}f(x,t)$. By solving $U(x,s)$ and using the
convolution theorem, we have its inverse Laplace transform:
$$u(x,t)=\mathscr{L}^{-1}U(x,s)=\mathscr{L}^{-1}\frac{1}{\prod^m_{i=1}(s^2-a_i^2P(\partial_x))}F(x,s)=G_m(\partial_x,t)*f(x,t).$$
By comparing (\ref{I2'}) with $u(x,t)=G_m(\partial_x,t)*f(x,t)$, we
have the expression of the abstract operators $G_m(\partial_x,t)$:
\begin{equation}\label{j6}
G_m(\partial_x,t)=\int^t_0\frac{(t-\tau)^{2m-3}}{(2m-3)!}
\sum^m_{j=1}\frac{a_j^{2m-2}}{\prod^m_{i=1 \atop
i\neq{j}}(a_j^2-a_i^2)}\frac{\sinh(\tau{a_j}P(\partial_x)^{1/2})}{a_jP(\partial_x)^{1/2}}\,d\tau.
\end{equation}
Applying (\ref{j6}) to (\ref{j5''}), thus Corollary 3.3 is proved.

\subsection{Analytic solutions of Cauchy problem}
\noindent

\textbf{Theorem 3.5.} Let
$\Delta_n=\frac{\partial^2}{\partial{x^2_1}}+\frac{\partial^2}{\partial{x^2_2}}+\cdots+\frac{\partial^2}{\partial{x^2_n}}$
be the n-dimensional Laplacian. If $n-2=2\nu+1,\;\nu\in\mathbb{N},\;a\in\mathbb{R}$,
then we have
\begin{equation}\label{36}
    \left\{\begin{array}{l@{\qquad}l}\displaystyle
    \left(\frac{\partial^2}{\partial{t^2}}-a^2\Delta_n\right)^mu=f(x,t),&x\in\mathbb{R}^n,\;t\in\mathbb{R}_+^1,\;\forall{m}\geq1,\\\displaystyle
    \left.\frac{\partial^ju}{\partial{t^j}}\right|_{t=0}=\varphi_j(x),&j=0,1,2,\ldots,2m-1.
    \end{array}\right.
\end{equation}
\begin{eqnarray}\label{36'}
u(x,t)
    &=& \int^t_0d\tau\int^{t-\tau}_0d\tau'\frac{\left[(t-\tau)^2-\tau'^2\right]^{m-2}\tau'^2}{(2m-2)!!\,(2m-4)!!}\nonumber\\
    & & \times\underbrace{\int^{\tau'}_0\tau'd\tau'\cdots}_\nu\int^{\tau'}_0\frac{(a^2\Delta_n)^\nu}{S'_n}
  \int_{S'_n}f(\xi',\tau)\,dS'_n\,\tau'd\tau' \nonumber\\
    & & +\,\frac{1}{(2m-2)!!}\sum^{\nu-1}_{r=0}\int^t_0\frac{(t-\tau)^{2m+2r-1}}{(2m+2r-1)!!\,(2r)!!}(a^2\Delta_n)^rf(x,\tau)\,d\tau \nonumber\\
    & & +\,\sum^{m-1}_{k=0}(-1)^k{m\choose{k}}(a^2\Delta_n)^{k+\nu}\sum^{2m-1-2k}_{j=0}\frac{\partial^{2m-1-2k-j}}{\partial{t}^{2m-1-2k-j}}
    \int^t_0d\tau\frac{(t^2-\tau^2)^{m-2}\tau^2}{(2m-2)!!\,(2m-4)!!} \nonumber\\
    & & \times\underbrace{\int^\tau_0\tau{d}\tau\cdots}_\nu\int^\tau_0\frac{1}{S_n}\int_{S_n}\varphi_j(\xi)\,dS_n\,\tau{d}\tau
    +\sum^{m-1}_{k=0}(-1)^k{m\choose{k}}\sum^{2m-1-2k}_{j=0}\nonumber\\
    & & \times\sum^{\nu-1}_{i=0}{m-1+i\choose{i}}\frac{t^{2k+2i+j}}{(2k+2i+j)!}(a^2\Delta_n)^{k+i}\varphi_j(x).
\end{eqnarray}
Where $n-2=2\nu+1$, $S'_n=2(2\pi)^{\nu+1}(a\tau')^{n-1}$,
$S_n=2(2\pi)^{\nu+1}(a\tau)^{n-1}$, and $\xi'\in\mathbb{R}_n$ is the
integral variable. The integral is on the hypersphere
$(\xi'_1-x_1)^2+(\xi'_2-x_2)^2+\cdots+(\xi'_n-x_n)^2=(a\tau')^2$,
and $dS'_n$ is its surface element. $\xi\in\mathbb{R}_n$ is the
integral variable on the hypersphere
$(\xi_1-x_1)^2+(\xi_2-x_2)^2+\cdots+(\xi_n-x_n)^2=(a\tau)^2$, and
$dS_n$ is its surface element.

\textbf{Proof.} In Theorem 3.2, let $P(\partial_x)=a^2\Delta_n$, applying (\ref{37}) to (\ref{32'}) we have
\begin{eqnarray*}
u(x,t)
    &=& \int^t_0\int^{t-\tau}_0\frac{\left[(t-\tau)^2-\tau'^2\right]^{m-2}}{(2m-2)!!\,(2m-4)!!}\nonumber\\
    & & \times\left(\tau'\underbrace{\int^{\tau'}_0\tau'd\tau'\cdots}_\nu\int^{\tau'}_0\frac{(a^2\Delta_n)^\nu}{S'_n}
  \int_{S'_n}f(\xi',\tau)\,dS'_n\,\tau'd\tau'\right)\tau'd\tau'd\tau \nonumber\\
    & & +\int^t_0\int^{t-\tau}_0\frac{\left[(t-\tau)^2-\tau'^2\right]^{m-2}}{(2m-2)!!\,(2m-4)!!}
    \sum^{\nu-1}_{r=0}\frac{\tau'^{2r+1}}{(2r+1)!}(a^2\Delta_n)^rf(x,\tau)\,\tau'd\tau'd\tau \nonumber\\
    & & +\,\sum^{m-1}_{k=0}(-1)^k{m\choose{k}}(a^2\Delta_n)^k\sum^{2m-1-2k}_{j=0}\frac{\partial^{2m-1-2k-j}}{\partial{t}^{2m-1-2k-j}}
    \int^t_0\frac{(t^2-\tau^2)^{m-2}\tau}{(2m-2)!!\,(2m-4)!!} \nonumber\\
    & & \times\left(\tau\underbrace{\int^\tau_0\tau{d}\tau\cdots}_\nu\int^\tau_0\frac{(a^2\Delta_n)^\nu}{S_n}\int_{S_n}\varphi_j(\xi)\,dS_n\,\tau{d}\tau\right)d\tau\nonumber\\
    & & +\,\sum^{m-1}_{k=0}(-1)^k{m\choose{k}}(a^2\Delta_n)^k\sum^{2m-1-2k}_{j=0}\frac{\partial^{2m-1-2k-j}}{\partial{t}^{2m-1-2k-j}}
    \int^t_0\frac{(t^2-\tau^2)^{m-2}\tau}{(2m-2)!!\,(2m-4)!!} \nonumber\\
    & & \times\sum^{\nu-1}_{i=0}\frac{\tau^{2i+1}}{(2i+1)!}(a^2\Delta_n)^i\varphi_j(x)d\tau.
\end{eqnarray*}
Where $n-2=2\nu+1$, $S'_n=2(2\pi)^{\nu+1}(a\tau')^{n-1}$,
$S_n=2(2\pi)^{\nu+1}(a\tau)^{n-1}$, and $\xi'\in\mathbb{R}_n$ is the
integral variable. The integral is on the hypersphere
$(\xi'_1-x_1)^2+(\xi'_2-x_2)^2+\cdots+(\xi'_n-x_n)^2=(a\tau')^2$,
and $dS'_n$ is its surface element. $\xi\in\mathbb{R}_n$ is the
integral variable on the hypersphere
$(\xi_1-x_1)^2+(\xi_2-x_2)^2+\cdots+(\xi_n-x_n)^2=(a\tau)^2$, and
$dS_n$ is its surface element.

By using (\ref{34}), where
$$\frac{\partial^{2m-1-2k-j}}{\partial{t}^{2m-1-2k-j}}\int^t_0\frac{(t^2-\tau^2)^{m-2}\tau}{(2m-2)!!\,(2m-4)!!}\frac{\tau^{2i+1}}{(2i+1)!}d\tau
={m-1+i\choose{i}}\frac{t^{2k+2i+j}}{(2k+2i+j)!}.$$
Similarly
$$\int^{t-\tau}_0\frac{\left[(t-\tau)^2-\tau'^2\right]^{m-2}}{(2m-2)!!\,(2m-4)!!}\frac{\tau'^{2r+1}}{(2r+1)!}\tau'd\tau'
=\frac{1}{(2m-2)!!}\frac{(t-\tau)^{2m+2r-1}}{(2m+2r-1)!!\,(2r)!!}.$$
So Theorem 3.5 is proved.

The Theorems 3.2 and Theorems 3.5 have been published by Reference \cite{bi11}.

\textbf{Theorem 3.6.} Let $a_1,a_2,\ldots,a_m$ be arbitrary real
roots different from each other, which satisfy
$\sum^m_{k=0}b_{2k}\chi^{2k}=\prod^m_{i=1}(\chi^2-a_i^2)$, and
$P(\partial_x,\partial_t)$ be the partial differential operators
defined by
$$P(\partial_x,\partial_t)=\sum^m_{k=0}b_{2k}\Delta_n^{m-k}\frac{\partial^{2k}}{\partial{t^{2k}}},
\quad{x}\in\mathbb{R}^n,\;t\in\mathbb{R}_+^1,\;\forall{m}\geq1.$$
Where $\Delta_n=\frac{\partial^2}{\partial{x^2_1}}+\frac{\partial^2}{\partial{x^2_2}}+\cdots+\frac{\partial^2}{\partial{x^2_n}}$ is the n-dimensional Laplacian. If
$n-2=2\nu+1,\;\nu=0,1,2,\cdots$, then we have
\begin{equation}\label{j7}
    \left\{\begin{array}{l@{\qquad}l}\displaystyle
    P(\partial_x,\partial_t)u=f(x,t),&x\in\mathbb{R}^n,\;t\in\mathbb{R}_+^1,\\\displaystyle
    \left.\frac{\partial^ru}{\partial{t^r}}\right|_{t=0}=\varphi_r(x),&r=0,1,2,\ldots,2m-1.
    \end{array}\right.
\end{equation}
\begin{eqnarray}\label{j7'}
  u(x,t) &=&
  \int^t_0d\tau\int^{t-\tau}_0d\tau'\frac{(t-\tau-\tau')^{2m-3}}{(2m-3)!}\sum^m_{j=1}\frac{a_j^{2m-2}}{\prod^m_{i=1 \atop i\neq{j}}(a_j^2-a_i^2)}\nonumber\\
  & &
  \times\left(\tau'\underbrace{\int^{\tau'}_0\tau'd\tau'\cdots}_\nu\int^{\tau'}_0\frac{(a_j^{2}\Delta_n)^{\nu}}{S'_{n,j}}\int_{S'_{n,j}}f(\xi',\tau)\,dS'_{n,j}\,\tau'd\tau'\right)\nonumber\\
  & & +\int^t_0\sum^m_{j=1}\frac{a_j^{2m-2}}{\prod^m_{i=1 \atop i\neq{j}}(a_j^2-a_i^2)}\sum^{\nu-1}_{l=0}\frac{a_j^{2l}\,(t-\tau)^{2l+2m-1}}{(2l+2m-1)!}\Delta_n^lf(x,\tau)d\tau\nonumber\\
  & & +\,\sum^m_{k=1}b_{2k}\Delta_n^{m-k}\sum^{2k-1}_{r=0}\frac{\partial^{2k-1-r}}{\partial t^{2k-1-r}}\int^t_0d\tau
  \frac{(t-\tau)^{2m-3}}{(2m-3)!}\sum^m_{j=1}\frac{a_j^{2m-2}}{\prod^m_{i=1 \atop i\neq{j}}(a_j^2-a_i^2)}\nonumber\\
  & &
  \times\left(\tau\underbrace{\int^{\tau}_0\tau{d\tau}\cdots}_\nu\int^{\tau}_0\frac{(a_j^{2}\Delta_n)^{\nu}}{S_{n,j}}\int_{S_{n,j}}\varphi_r(\xi)\,dS_{n,j}\,\tau{d\tau}\right)\nonumber\\
  & & +\,\sum^m_{k=1}b_{2k}\Delta_n^{m-k}\sum^{2k-1}_{r=0}\sum^m_{j=1}\frac{a_j^{2m-2}}{\prod^m_{i=1 \atop i\neq{j}}(a_j^2-a_i^2)}\sum^{\nu-1}_{l=0}\frac{a_j^{2l}\,t^{2l+2m+r-2k}}{(2l+2m+r-2k)!}\Delta_n^l\varphi_r(x).
\end{eqnarray}
Where $S'_{n,j}=2(2\pi)^{\nu+1}(a_j\tau')^{n-1}$,
$S_{n,j}=2(2\pi)^{\nu+1}(a_j\tau)^{n-1}$, and $\xi'\in\mathbb{R}_n$
is the integral variable. The integral is on the hypersphere
$(\xi'_1-x_1)^2+(\xi'_2-x_2)^2+\cdots+(\xi'_n-x_n)^2=(a_j\tau')^2$,
and $dS'_{n,j}$ is its surface element. $\xi\in\mathbb{R}_n$ is the
integral variable on the hypersphere
$(\xi_1-x_1)^2+(\xi_2-x_2)^2+\cdots+(\xi_n-x_n)^2=(a_j\tau)^2$, and
$dS_{n,j}$ is its surface element.

\textbf{Proof.} In Theorem 3.4, let $P(\partial_x)=\Delta_n$, applying (\ref{37}) to (\ref{j2'}) we have
\begin{eqnarray*}
  u(x,t) &=&
  \int^t_0d\tau\int^{t-\tau}_0d\tau'\frac{(t-\tau-\tau')^{2m-3}}{(2m-3)!}\sum^m_{j=1}\frac{a_j^{2m-2}}{\prod^m_{i=1 \atop i\neq{j}}(a_j^2-a_i^2)}\nonumber\\
  & &
  \times\left(\tau'\underbrace{\int^{\tau'}_0\tau'd\tau'\cdots}_\nu\int^{\tau'}_0\frac{(a_j^{2}\Delta_n)^{\nu}}{S'_{n,j}}\int_{S'_{n,j}}f(\xi',\tau)\,dS'_{n,j}\,\tau'd\tau'
  +\sum^{\nu-1}_{l=0}\frac{a_j^{2l}\tau'^{2l+1}}{(2l+1)!}\Delta_n^lf(x,\tau)\right)\nonumber\\
   & & +\,\sum^m_{k=1}b_{2k}\Delta_n^{m-k}\sum^{2k-1}_{r=0}\frac{\partial^{2k-1-r}}{\partial t^{2k-1-r}}\int^t_0d\tau
   \frac{(t-\tau)^{2m-3}}{(2m-3)!}\sum^m_{j=1}\frac{a_j^{2m-2}}{\prod^m_{i=1 \atop i\neq{j}}(a_j^2-a_i^2)}\nonumber\\
   & &
   \times\left(\tau\underbrace{\int^{\tau}_0\tau{d\tau}\cdots}_\nu\int^{\tau}_0\frac{(a_j^{2}\Delta_n)^{\nu}}{S_{n,j}}\int_{S_{n,j}}\varphi_r(\xi)\,dS_{n,j}\,\tau{d\tau}
   +\sum^{\nu-1}_{l=0}\frac{a_j^{2l}\tau^{2l+1}}{(2l+1)!}\Delta_n^l\varphi_r(x)\right).
\end{eqnarray*}
Where $S'_{n,j}=2(2\pi)^{\nu+1}(a_j\tau')^{n-1}$,
$S_{n,j}=2(2\pi)^{\nu+1}(a_j\tau)^{n-1}$, and $\xi'\in\mathbb{R}_n$
is the integral variable. The integral is on the hypersphere
$(\xi'_1-x_1)^2+(\xi'_2-x_2)^2+\cdots+(\xi'_n-x_n)^2=(a_j\tau')^2$,
and $dS'_{n,j}$ is its surface element. $\xi\in\mathbb{R}_n$ is the
integral variable on the hypersphere
$(\xi_1-x_1)^2+(\xi_2-x_2)^2+\cdots+(\xi_n-x_n)^2=(a_j\tau)^2$, and
$dS_{n,j}$ is its surface element.
Then Theorem 3.6 is proved.

Similarly, we can easily obtain explicit solutions of the Cauchy problem of more complex partial differential equations, thus we have established the
general theory of initial value problems for linear higher-order partial differential equations.

\subsection{Solvability for initial-boundary value problem}
\noindent

Clearly, we can attach proper boundary conditions to the initial
value problem introduced by Theorem 3.2, Theorem 3.3 and Theorem 3.4, which makes
the given functions $f(x,t),\,\varphi_j(x)$ in Eq (\ref{32}), (\ref{j1}) and
(\ref{j2}) become functions with boundary conditions.
For the initial-boundary value problem, the partial differential operators $P(\partial_x)$
must have the characteristic function related to boundary
conditions, in order to the known functions
$f(x,\tau),\,\varphi_r(x)$ in (\ref{32'}), (\ref{j1'}) and
(\ref{j2'}) can be expanded in infinite series of characteristic
function of $P(\partial_x)$. So the $P(\partial_x)$ in Theorem 3.2,
Theorem 3.3 and Theorem 3.4 can be variable-coefficient partial
differential operators. In order to solve the corresponding
initial-boundary value problem, we need to solve the characteristic
value problem of $P(\partial_x)$ under given boundary conditions to
determine a set of orthogonal functions. For instance, in (\ref{32'}), (\ref{j1'}) and (\ref{j2'}), if
$f(x,\tau),\varphi_r(x)\in{L^2}(\Omega)$, and $P(\partial_x)$ is the
second-order linear self-adjoint elliptic operators, namely
\begin{equation}\label{38'}
  P(\partial_x)u=\sum^n_{i,j=1}\frac{\partial}{\partial{x}_j}\left(a_{ij}(x)\frac{\partial{u}}{\partial{x}_i}\right)+c(x)u,
\quad{x}\in\Omega\subset\mathbb{R}^n,
\end{equation}
then the boundary conditions can be added for the definite solution
problems (\ref{32}), (\ref{j1}) and (\ref{j2}):
$\overline{B}u|_{\partial\Omega}=0$, representing
$u|_{\partial\Omega}=0$ or
\[\left[\sum^n_{i,j=1}a_{ij}(x)\frac{\partial{u}}{\partial{x}_j}\cos\langle\mathbf{a},x_i\rangle+b(x)u\right]_{\partial\Omega}=0.\]
Where $\mathbf{a}$ is the unit outward normal of $\partial\Omega$.
Thus this kind of initial boundary value problems boils down to
solving the characteristic value problem of the first boundary value
problem of second-order linear self-adjoint elliptic operators:
\begin{equation}\label{38}
    \left\{\begin{array}{l@{\qquad}l}\displaystyle
    \sum^n_{i,j=1}\frac{\partial}{\partial{x}_j}\left(a_{ij}(x)\frac{\partial{u}}{\partial{x}_i}\right)+c(x)u=
    -\lambda{u},\quad{x}\in\Omega\subset\mathbb{R}^n,\\\displaystyle
    \overline{B}u|_{\partial\Omega}=0.
    \end{array}\right.
\end{equation}

For the characteristic value problem (\ref{38}), by Reference \cite{evans} we have the following results:

Let $\Omega\subset\mathbb{R}^n$ be a bounded open domain, and
$\partial\Omega$ be smooth. Let $a_{ij}=a_{ji}$, and there exists
$\theta>0$ such that
$$\sum^n_{i,j=1}a_{ij}(x)\xi_i\xi_j\geq\theta|\xi|^2,\quad{x\in\Omega}.$$
And let
$a_{ij}\in{C}^1(\overline{\Omega}),\;c(x)\in{C}(\overline{\Omega}),\;b(x)\in{C}(\partial\Omega)$,
then (\ref{38}) has the following countable characteristic values:
\[0\leq\lambda_1\leq\lambda_2\leq\cdots\leq\lambda_\nu\leq\cdots,\quad\lim_{\nu\rightarrow\infty}\lambda_\nu=\infty\]
(If $(a_{ij})=I$ is a unit matrix, then $\lambda_1=0$ when
$b(x)=c(x)=0$. When $b(x)\geq0,\,c(x)\geq0$ and one of them does not
identically equal to zero, $\lambda_1>0$) and the corresponding
characteristic functions $e_1(x),e_2(x),\cdots,e_\nu(x),\cdots,$
satisfy
\begin{equation}\label{39}
  \sum^n_{i,j=1}\frac{\partial}{\partial{x}_j}\left(a_{ij}(x)\frac{\partial{e_i}}{\partial{x}_i}\right)+c(x)e_i=-\lambda_i{e_i},\quad
(e_i,\,e_j)=\delta_{ij}
\end{equation}
and $\{e_j(x)\}^\infty_{j=1}$ is complete in $L^2(\Omega)$, that is
for an arbitrary $f(x)\in{L^2}(\Omega)$, there exists $c_j$ such that
\[\lim_{\nu\rightarrow\infty}\|f-\sum^\nu_{i=1}c_ie_i\|_{L^2(\Omega)}=0.\]
Particularly, let $\Omega\in\mathbb{R}^n$ be a bounded smooth domain, then for the characteristic value problem of the Laplace operators
\begin{equation}\label{ws}
    \left\{\begin{array}{l@{\qquad}l}\displaystyle
    \Delta_nu=-\lambda{u},&x\in\Omega,\\\displaystyle
    u|_{\partial\Omega}=0,
    \end{array}\right.
\end{equation}
 an orthogonal system of $H_0^1(\Omega)$ is composed of its solutions $\{e_j(x)\}^\infty_{j=1}$ (See \cite{ws}). Here $H_0^1(\Omega)=W_0^{1,2}(\Omega),\;W_0^{m,p}(\Omega)
 =\{u\in{W}^{m,p}(\Omega)\mid\mbox{exists}\;u_k\in{C}_0^\infty(\Omega),\,\mbox{making}\;\|u_k-u\|_{m,p,\Omega}\rightarrow0(k\rightarrow\infty)\}$. $W^{m,p}(\Omega)$ is the $m$-order Sobolev space on $\Omega$.

Therefore, for (\ref{32'}), (\ref{j1'}) and (\ref{j2'}), by
$f(x,\tau),\varphi_r(x)\in{L^2}(\Omega)$ we also have
\[\lim_{\nu\rightarrow\infty}\|f(x,\tau)-\sum^\nu_{i=1}c_i(\tau)e_i(x)\|_{L^2(\Omega)}=0.\]
\[\lim_{\nu\rightarrow\infty}\|\varphi_r(x)-\sum^\nu_{i=1}c_ie_i(x)\|_{L^2(\Omega)}=0.\]
Clearly, based on (\ref{38'}) and (\ref{39}), the abstract operators
$$e^{t{a_j}P(\partial_x)},\quad\frac{\sinh(ta_jP(\partial_x)^{1/2})}{a_jP(\partial_x)^{1/2}}\quad\mbox{and}\quad\cosh(ta_jP(\partial_x)^{1/2})
=\frac{\partial}{\partial{t}}\,
\frac{\sinh(ta_jP(\partial_x)^{1/2})}{a_jP(\partial_x)^{1/2}}$$ becomes the
continuous operators on the Hilbert space, which are also the linear continuous mapping $f(t,P(\partial_x))$ :
$L^2((0,T]\times\Omega)\rightarrow{L^2((0,T]\times\Omega)}$, and they
acts on the characteristic function $e_i(x)$, we have
\begin{equation}\label{40}
  e^{\tau{a_j}P(\partial_x)}e_i(x)=e^{-\tau{a_j}\lambda_i}e_i(x),\quad\frac{\sinh(\tau{a_j}P(\partial_x)^{1/2})}{a_jP(\partial_x)^{1/2}}e_i(x)=
  \frac{\sin(a_j\sqrt{\lambda_i}\,\tau)}{a_j\sqrt{\lambda_i}}e_i(x),
\end{equation}
where $x\in\Omega\subset\mathbb{R}^n,\;j=1,2,\cdots,m$.
So we have
\[\lim_{\nu\rightarrow\infty}\|e^{\tau'{a_j}P(\partial_x)}f(x,\tau)-\sum^\nu_{i=1}(f,e_i)e^{-\tau'{a_j}\lambda_i}e_i(x)\|_{L^2(\Omega)}=0.\]
\[\lim_{\nu\rightarrow\infty}\|e^{\tau{a_j}P(\partial_x)}\varphi_r(x)-\sum^\nu_{i=1}(\varphi_r,e_i)e^{-\tau{a_j}\lambda_i}e_i(x)\|_{L^2(\Omega)}=0.\]
\[\lim_{\nu\rightarrow\infty}\left\|\frac{\sinh(\tau'a_jP(\partial_x)^{1/2})}{a_jP(\partial_x)^{1/2}}f(x,\tau)-
\sum^\nu_{i=1}(f,e_i)\frac{\sin(a_j\sqrt{\lambda_i}\,\tau')}{a_j\sqrt{\lambda_i}}e_i(x)\right\|_{L^2(\Omega)}=0.\]
\[\lim_{\nu\rightarrow\infty}\left\|\frac{\sinh(\tau{a_j}P(\partial_x)^{1/2})}{a_jP(\partial_x)^{1/2}}\varphi_r(x)
-\sum^\nu_{i=1}(\varphi_r,e_i)\frac{\sin(a_j\sqrt{\lambda_i}\,\tau)}{a_j\sqrt{\lambda_i}}e_i(x)\right\|_{L^2(\Omega)}=0.\]
Thus when (\ref{32}), (\ref{j1}) and (\ref{j2}) satisfy the boundary
condition $\overline{B}u|_{\partial\Omega}=0$, and if the orthogonal
complete series $\{e_j(x)\}^\infty_{j=1}$ in the Hilbert space are known, we get the
classical solution of the corresponding initial-boundary value problem.


\begin{thebibliography}{8}
\bibitem{chen} S.X. Chen, 2006. Pseudodifferential Operators, (Second Edition). Beijing: Higher Education Press. China. pp. 8-25
\bibitem{bi97} G.Q. Bi, Applications of abstract operators in partial differential equation({\romannumeral1}), \emph{Pure and Applied Mathematics}. 1997, \textbf{13}(1): 7-14
\bibitem{bi99} G.Q. Bi, Applications of abstract operators in partial differential equation({\romannumeral2}), \emph{Chin. Quart. J. of Math.} 1999, \textbf{14}(3): 80-87
\bibitem{bi01} G.Q. Bi, Operator methods in high order partial differential equation, \emph{Chin. Quart. J. of Math.} 2001, \textbf{16}(1): 88-101
\bibitem{bi11} G.Q. Bi and Y.K. Bi, Abstract operators and higher-order linear partial differential equation, \emph{Chin. Quart. J. of Math.} 2011, \textbf{26}(4): 511-515
\bibitem{bigy} G.Q. Bi and Y.K. Bi, New properties of Fourier series and Riemann Zeta function, 2010. arXiv:1008.5046
\bibitem{evans} L.C. Evans, 1998. Partial Differential Equations. Graduate Studies in Mathematics Vol.19. Rhode Island: American Mathematical Society.
\bibitem{ws} S. Wang, 2009. Sobolev Space and Partial Differential Equations Introdunction. Beijing: Science Press. China. pp. 155-157
\end{thebibliography}
\end{document}